\documentclass[final]{article} 
\usepackage{amsthm,amsmath,amsfonts,amssymb,xspace,paralist,pictexwd}
\usepackage[latin1]{inputenc}
\usepackage{textcomp}
\usepackage{mathrsfs}
\usepackage{fixme,proof}  \usepackage{mathptmx}

\usepackage[scaled=0.92]{helvet}

\newtheorem{theorem}{Theorem}
\newtheorem{lemma}[theorem]{Lemma}
\newtheorem{proposition}[theorem]{Proposition}
\newtheorem{corollary}[theorem]{Corollary}

\theoremstyle{definition}
\newtheorem{definition}[theorem]{Definition}
\newtheorem{notation}[theorem]{Notation}

\theoremstyle{remark}
\newtheorem*{remark}{Remark}
\newtheorem*{example}{Example}

\let\imm\ensuremath
\newcommand\limp{\to}
\let\impl\limp
\newcommand{\rimpl}{\leftarrow}

\newcommand\I{\imm{\mathfrak{I}}}
\newcommand\J{\imm{\mathfrak{J}}}
\def\boeta{{\bf\eta}}
\def\one{{{\bf 1}}}

\newcommand{\g}{$\mathbf{G}$}
\newcommand\LC{\imm{\mathbf{LC}}}

\newcommand\Gd{{\gdl V}}
\newcommand\Go{{\mathbf {G}^0}}
\newcommand\Gn{{\gdl n}}
\newcommand\GR{{\gdl\R}}
\newcommand\R{\mathbb{R}}
\newcommand\dom{\mathop{\mathrm{dom}}}
\def\HU{\mathrm{HU}}
\def\HB{\mathrm{HB}}
\let\bbR\R

\let\vup\up
\def\dn{{\downarrow}}
\let\vdn\dn

\def\card#1{\lvert #1 \rvert}
\def\eval#1{\lvert #1 \rvert}

\def\T{{\bf\sf T}}

\newcommand{\modelsg}{\mathrel{\Gd \models}}

\def\N{{\mathbb N}}
\def\bbZ{{\mathbb Z}}
\def\bbD{\mathbb{D}}
\let\bbN\N
\let\nmodels\nvDash
\def\proves{\mathrel{\vdash}}
\def\nproves{\mathrel{\nvdash}}
\def\entails{\mathrel{\models}}
\def\nentails{\mathrel{\nmodels}}
\def\dequiv{\Leftrightarrow}

\def\strc#1#2{\left<#1,#2\right>}
\def\Feq{{\cF/\mathnormal{\equiv}}}
\def\lF{\strc\Feq\le}

\newcommand\tup[1]{\overline{#1}}
\def\Q{{\ensuremath{\mathord{\mathsf{Q}}}}}
\newcommand{\LL}{{\mathscr L}}

\newcommand\qe[1]{\exists #1}
\newcommand\qa[1]{\forall #1}
\newcommand\suchthat{\mathrel{:}}

\newcommand{\gdl}[1]{\mathbf {G}_{#1}}
\newcommand\subval{\mathrm{Val}}
\newcommand{\bbQ}{\ensuremath{\mathbb{Q}}\xspace}

\newcommand{\IPL}{\mbox{IPL}}
\newcommand{\IL}{\mbox{IL}}

\renewcommand{\H}{\mathbf{H}}
\newcommand{\Ho}{{\mathbf{H}_0}}

\newcommand{\cF}{\mathcal{F}}

\newcommand{\cT}{\mathcal{T}}
\newcommand{\cA}{\mathcal{A}}

\def\axlin{{\ensuremath{\hbox{\scshape lin}}}}
\def\axiso{{\ensuremath{\hbox{\scshape iso}}}}
\def\axqs{{\ensuremath{\hbox{\scshape qs}}}}
\def\axqs{{\ensuremath{\hbox{\scshape qs}}}}

\def\axlin{{\ensuremath{\hbox{\scshape lin}}}}
\def\axiso{{\ensuremath{\hbox{\scshape iso}_0}}}

\def\axfinn{{\ensuremath{\hbox{\scshape fin}(n)}}}

\def\one{\mathbf{1}}
\def\zero{\mathbf{0}}
\def\los#1{\strc{#1}\le}

\def\lR{\los{[0,1]}}
\def\alepho{{\aleph_0}}
\def\lee{\preceq}

\def\lAnd{\bigwedge}

\let\Proves\Rightarrow
\let\nProves\nRightarrow

\def\tik#1:{%
  \putrule from #1 -2 to #1 +2}
\def\Tik#1,#2:{%
  \put {$#2$} [b] at #1 5
  \putrule from #1 -2 to #1 +4}
\def\tiku#1,#2:{%
  \put {$#2$} [b] at #1 24
  \putrule from #1 18 to #1 22}
\def\tikd#1,#2:{%
  \put {$\smash{#2}$} at #1 -6
  \putrule from #1 -2 to #1 2}
\def\tikdd#1,#2:{%
  \put {$\smash{#2}$} at #1 4
  \putrule from #1 -2 to #1 2}

\title{First-Order Gödel Logics}

\author{Matthias Baaz\thanks{Research supported by \textsc{FWF} grant
\textsc{P15477--MAT}}\\ 
Technische Universität Wien\\ 
1040 Vienna, Austria\\
baaz@logic.at \and
Norbert Preining%
\thanks{Research supported by \textsc{EC-MC 008054}
  and \textsc{FWF} grant \textsc{P16539-N04}}\\
Università di Siena\\
53100 Siena, Italy\\
preining@logic.at \and
Richard Zach\thanks{Research supported by the Natural Sciences and Engineering
  Research Council of Canada}\\
University of Calgary\\
Calgary, \textsc{AB T2N 1N4},
Canada\\ rzach@ucalgary.ca}

\begin{document}

\maketitle

\begin{abstract} 
  First-order Gödel logics are a family of infinite-valued logics where the
  sets of truth values~$V$ are closed subsets of $[0, 1]$ containing both $0$
  and $1$.  Different such sets~$V$ in general determine different Gödel
  logics~$\gdl{V}$ (sets of those formulas which evaluate to~$1$ in every
  interpretation into~$V$).  It is shown that $\gdl{V}$ is axiomatizable iff
  $V$ is finite, $V$ is uncountable with $0$ isolated in $V$, or every
  neighborhood of $0$ in $V$ is uncountable.  Complete axiomatizations for
  each of these cases are given.  The r.e. prenex, negation-free, and
  existential fragments of all first-order Gödel logics are also
  characterized.
\end{abstract}
\tableofcontents

\section{Introduction}

\subsection{Motivation}

The logics we investigate in this paper, first-order Gödel logics, can be
characterized in a rough-and-ready way as follows:  The language is a standard
first-order language.  The logics are many-valued, and the sets of truth
values considered are closed subsets of $[0, 1]$ which contain both~$0$
and~$1$. $1$ is the ``designated value,'' i.e., a formula is valid if
it receives the 
value~$1$ in every interpretation.  The truth functions of conjunction and
disjunction are minimum and maximum, respectively, and quantifiers are defined
by infimum and supremum over subsets of the set of truth values.  The
characteristic operator of Gödel logics, the Gödel conditional, is defined by
$a \limp b = 1$ if $a \le b$ and $= b$ if $a > b$.  Because the truth values
are ordered (indeed, in many cases, densely ordered), the semantics of Gödel
logics is suitable for formalization of \emph{comparisons}.  It is
related in this respect to a more widely known many-valued logic, \L
ukasiewicz (or ``fuzzy'') logic---yet the truth function of the \L ukasiewicz
conditional is defined not just using comparison, but also addition.  In
contrast to \L ukasiewicz logic, which might be considered a logic of
\emph{absolute} or \emph{metric comparison}, Gödel logics are logics of
\emph{relative comparison}.  This alone makes Gödel logics an interesting
subject for logical investigations.

There are other reasons why the study of Gödel logics is important.  As noted,
Gödel logics are related to other many-valued logics of recognized importance.
Indeed, Gödel logic is one of the three basic
$t$-norm based logics which have received increasing attention in the last 15
or so years \cite{hajek} (the others are \L ukasiewicz and product logic).
Yet Gödel logic is also closely related to intuitionistic logic: it is the
logic of linearly-ordered Heyting algebras.  In the propositional case,
infinite-valued Gödel logic can be axiomatized by the intuitionistic
propositional calculus extended by the axiom schema $(A \limp B) \lor (B \limp
A)$.  This connection extends also to Kripke semantics for intuitionistic
logic: Gödel logics can also be characterized as logics of (classes of)
linearly ordered and countable intuitionistic Kripke structures with
constant domains \cite{BeckPrei05Kripke}.

One of the surprising facts about Gödel logics is that whereas there is
only one infinite-valued propositional Gödel logic, there are infinitely many
different infinite-valued first-order Gödel logics depending on the choice of
the set of truth values.  This is also the case when one considers the
propositional consequence relation, and likewise when the language is extended
to include quantification over propositions.  For both quantified
propositional and first-order Gödel logics, different sets of truth values
with different order-theoretic properties result in different sets of valid
formulas.  Hence it is necessary to consider truth value sets other than the
standard unit interval.

In the light of the result of Scarpellini \cite{scarpellini} on
non-axiomatizability of infinite-valued first-order \L ukasiewicz logic
which can be extended to almost all linearly ordered infinite-valued 
logics, it is also surprising that some infinite-valued Gödel logics are 
recursively enumerable. Our
main aim in this paper is to characterize those sets of truth values which
give rise to axiomatizable Gödel logics, and those whose sets of validities
are not r.e.  We show that a set $V$ of truth values determines an
axiomatizable first-order Gödel logic if, and only if, $V$ is finite, $V$ is
uncountable and $0$ is isolated, or every neighborhood of $0$ in $V$ is
uncountable.  These cases also determine different sets of validities: the
finite-valued Gödel logics $\Gn$, the logic $\Go$, and the ``standard''
infinite-valued Gödel logic~$\GR$ (based on the truth value set $[0, 1]$).

\subsection{History of Gödel logics}

Gödel logics are one of the oldest families of many-valued logics.
Propositional finite-valued Gödel logics were introduced by Gödel in
\cite{goedel33} to show that intuitionistic logic does not have a
characteristic finite matrix. They provide the first examples of intermediate
logics (intermediate, that is, in strength between classical and
intuitionistic logics).  Dummett \cite{dummett} was the first to study
infinite valued propositional Gödel logics, axiomatizing the set of
tautologies over infinite truth-value sets by intuitionistic logic extended by
the linearity axiom $(A \limp B) \lor (B \limp A)$. Hence, infinite-valued
propositional Gödel logic is also sometimes called Gödel-Dummett logic or
Dummett's \textsc{LC}. In terms of Kripke semantics, the characteristic
linearity axiom picks out those accessibility relations which are linear
orders.

Standard first-order Gödel logic~$\GR$---the one based on the full interval
$[0,1]$---has been discovered and studied by several people independently.
Alfred Horn was probably the first: He discussed this logic under the name
\emph{logic with truth values in a linearly ordered Heyting algebra}
\cite{horn}, and gave an axiomatization and the first completeness proof.
Takeuti and Titani~\cite{TT} called $\GR$ \emph{intuitionistic fuzzy logic},
and also gave an axiomatization for which they proved the completeness. This
system incorporates the density rule
\[ 
\infer{\Gamma\proves{A\lor (C\limp B)}}{\Gamma\proves{A\lor (C\limp p) \lor
    (p\limp B)}}
\]
(where~$p$ is any propositional variable not occurring in the lower sequent.)
The rule is redundant for an axiomatization of~$\GR$, as was shown by Takano
\cite{takano}, who gave a streamlined completeness proof of Takeuti-Titani's
system without the rule. (A syntactical proof of the elimination of the
density rule was later given in~\cite{BaazZach00CSL}.  Other proof-theoretic
investigations of Gödel logics can be found in \cite{BC02} and \cite{BFC03}.)
The density rule is nevertheless interesting: It forces the truth value set to
be dense in itself (in the sense that, if the truth value set isn't dense in
itself, the rule does not preserve validity).  This contrasts with the
expressive power of formulas: no formula is valid only for truth value sets
which are dense in themselves.

First-order Gödel logics other than $\GR$ were first considered in
\cite{BaazLeitschZach:96}, where it was shown that $\gdl{\vdn}$, based on the
truth value set $V_\vdn = \{1/k : k \in \N\} \cup \{0\}$ is not r.e.  Hájek
  \cite{Hajek05} has recently improved this result, and showed that not only
  is the set of validities not r.e., it is not even arithmetical.  Hájek
  also showed that the Gödel logic~$\gdl{\vup}$ based on $V_\vup =  \{1-1/k : k \in \N\}
    \cup \{1\}$ is $\Pi_2$-complete. Results preliminary to the results of the
    present paper were reported in \cite{BPZ03,Prei02LPAR,Prei03PHD}. 

\subsection{Overview of the results}

We begin with a preliminary discussion of the syntax and semantics of Gödel
logics, including a discussion of some of the more interesting special cases
of first-order Gödel logics and their relationships
(Section~\ref{sec:preliminaries}).  In Section~\ref{sec:order}, we
present some relevant results regarding the topology of truth-value
sets. 

The main results of the paper are contained in
Sections~\ref{sec:countable}--\ref{sec:finite}.  We provide a complete
classification of the axiomatizability of first order Gödel logics. The main
results are, that a logic based on a truth value set~$V$ is axiomatizable if
and only if
\begin{enumerate}
\item $V$ is finite (Section~\ref{sec:finite}), or 
\item $V$ is uncountable and $0$ is contained in the perfect kernel
  (Section~\ref{ssec:perfect}), or
\item $V$ is uncountable and $0$ is isolated (Section~\ref{ssec:isolated}).
\end{enumerate}
In all other cases, i.e., logics with countable truth value set
(Section~\ref{sec:countable}) and those where there is a countable
neighborhood of~$0$ and $0$ is not isolated (Section~\ref{sec:uncount-nonax}),
the respective logics are not r.e.

In Section~\ref{sec:fragments}, we investigate the complexity of fragments of
first-order Gödel logic, specifically, the prenex fragments
(Section~\ref{sec:prenex}), the $\bot$-free fragments
(Section~\ref{sec:botfree}), and the existential ($\forall$-free) fragments
(Section~\ref{sec:efrag}).  We show that the prenex fragment of a Gödel logic
is axiomatizable if and only if the truth value set is finite or uncountable.
This means that there are truth-value sets where the prenex fragment of the
corresponding logic is r.e. even though the full logic is not. Moreover, there
all axiomatizable prenex fragments coincide.  This is also the case for
$\bot$-free and existential fragments, but in these cases only those truth
value sets determine r.e. $\bot$-free and existential fragments for which also
the full logic is r.e., viz., truth value sets which are finite, uncountable
with $0$ isolated, and those where every neighborhood of $0$ is uncountable.

\section{Preliminaries}\label{sec:preliminaries}

\subsection{Syntax and Semantics}

In the following we fix a standard first-order language~$\LL$ with finitely or
countably many predicate symbols~$P$ and finitely or countably many function
symbols~$f$ for every finite arity~$k$.  In addition to the two quantifiers
$\forall $ and $\exists$ we use the connectives $\vee$, $\wedge$, $\to$ and
the constant $\bot$ (for `false'); other connectives are introduced as
abbreviations, in particular we let $\lnot A \equiv (A \to \bot)$.

Gödel logics are usually defined using the single truth value set~$[0,1]$. For
propositional logic the choice of any infinite subset of $[0,1]$ leads to the
same propositional logic (set of tautologies). In the first order case, where
quantifiers will be interpreted as infima and suprema, a closed subset of
$[0,1]$ is necessary.

\begin{definition}[Gödel set]
  A \emph{Gödel set} is a closed set~$V \subseteq [0,1]$
  which contains $0$ and~$1$.  
\end{definition}

The semantics of Gödel logics, with respect to a fixed Gödel set as truth
value set and a fixed language~$\LL$ of predicate logic, is defined using the
extended language $\LL^U$, where $U$ is the universe of the
interpretation~$\I$.  $\LL^U$ is $\LL$ extended with constant symbols for each
element of $U$.

\begin{definition}[Semantics of Gödel logic]\label{def:gsemantik}
  Fix a Gödel set $V$.
  An \emph{interpretation $\I$} into $V$ consists of 
  \begin{enumerate}
  \item  a nonempty set $U = U^\I$, the
    `universe' of $\I$, 
  \item  for each $k$-ary predicate symbol $P$, 
    a function $P^\I : U^k \to V$,
  \item for each $k$-ary function symbol $f$, a function
    $f^\I: U^k \to U$.
  \item for each variable~$v$, a value $v^\I \in U$.
  \end{enumerate}
  
  Given an interpretation $\I$, we can naturally define a value $t^\I$ for any
  term $t$ and a truth value $\I(A)$ for any formula $A$ of $\LL^U$. For a
  terms $t = f(u_1, \ldots, u_k)$ we define $\I(t)=f^\I(u_1^\I, \ldots,
  u_k^\I)$. For atomic formulas $A \equiv P(t_1,\dots, t_n)$, we define $\I(A)
  = P^\I(t_1^\I, \ldots, t_n^\I)$. For composite formulas $A$ we define
  $\I(A)$ by:
  \begin{align}
    \I(\bot) &= 0\\
    \I(A \land B) &= \min(\I(A),\I(B))\\
    \I(A \lor B) &= \max(\I(A),\I(B))\\
    \I(A \limp B) &= \begin{cases} 1 & \I(A)\le \I(B)\\
                                  \I(B) & \text{ otherwise}
                     \end{cases}\\
    \I(\qa x\, A(x)) &= \inf \{\I(A(u)) \suchthat u \in U\}\\
    \I(\qe x\, A(x)) &= \sup \{\I(A(u)) \suchthat u \in U\}
  \end{align}
  (Here we use the fact that every Gödel sets $V$ is a \emph{closed} subset of
  $[0,1]$ in order to be able to interpret $\forall$ and $\exists$ as $\inf$
  and $\sup$ in~$V$.)

  If $\I(A) = 1$, we say that $\I$ \emph{satisfies} $A$, and write $\I
  \entails A$.
\end{definition}

\begin{definition}[Gödel logics based on $V$]\label{def:goedellogics}
  For a Gödel set $V$ we define the {\em first order Gödel logic $\gdl V$} as
  the set of all formulas of $\LL$ such that $\I \entails A$ for all
  $V$-interpretations~$\I$.
\end{definition}

It should be noted that for Gödel logics with~$0$ isolated, the notion of 
\emph{satisfiability}
for sets of formulas is not particularly interesting, since a set of
formulas~$\Gamma$ is satisfiable (in the sense that there is an $\I$ so that
$\I \entails A$ for all $A \in \Gamma$) iff it is satisfiable classically.
\fxnote{$0$ nicht isoliert: $(\forall xP(x)\limp\bot)\land\forall x¬¬P(x)$ Gegenbsp!}
For this reason, we take \emph{entailment} to
be the fundamental model-theoretic notion.

\begin{definition}
  If $\Gamma$ is a set of formulas (possibly infinite), we say that $\Gamma$
  entails $A$ in $\gdl V$, $\Gamma \entails_V A$ iff for all $\I$ into $V$,
  \[ 
  \inf \{\I(B) : B \in \Gamma\} \le \I(A);
  \]
  and $\Gamma$ \emph{$1$-entails} $A$ in $\gdl V$, $\Gamma \Vdash_V A$, iff,
  for all~$\I$ into $V$, whenever $\I(B) = 1$ for all $B \in \Gamma$, then
  $\I(A) = 1$.
\end{definition}

\begin{notation}
  We will write $\Gamma \entails A$ instead of $\Gamma\entails_V A$ in
  case it is obvious which truth value set $V$ is meant. We will
  sometimes write $\Gamma \entails \Delta \in \gdl V$, by which we mean 
  that $\Gamma \entails_V \Delta$. The notation $\gdl V \entails A$
  stands for $\emptyset \entails_V A$, or $A \in \gdl V$. 
\end{notation}

Whether or not a formula $A$ evaluates to~1 under an interpretation~$\I$
depends only on the \emph{relative ordering} of the truth values of the atomic
formulas (in $\LL^\I$), and not directly on the set~$V$ or on the
\emph{values} of the atomic formulas. If $V \subseteq W$ are both Gödel sets,
and $\I$ is an interpretation into $V$, then $\I$ can be seen also as a
interpretation into $W$, and the values $\I(A)$, computed recursively using
(1)--(6), do not depend on whether we view $\I$ as a $V$-interpretation or a
$W$-interpretation.  Consequently, if $V \subseteq W$, there are more
interpretations into $W$ than into $V$.  Hence, if $\Gamma \entails_W A$ then
also $\Gamma \entails_V A$ and $\gdl{W} \subseteq \gdl{V}$.

This can be generalized to embeddings between Gödel sets other than inclusion.
First, we make precise which formulas are involved in the computation of the
truth-value of a formula~$A$ in an interpretation~$\I$:

\begin{definition}\label{def:subformula}
  The only subformula of an atomic formula $P$ in $\LL^U$ is  $P$
  itself. The subformulas of $A \star B$ for
  $\star \in \{\limp,\land,\lor\}$ are the subformulas of  $A$ and of
  $B$, together with $A \star B$ itself.  The subformulas of 
  $\qa x\,A(x)$ and $\qe x\,A(x)$ with respect to a universe $U$ are all
  subformulas of all $A(u)$ for $u\in U$, together with  $\qa x\,A(x)$
  (or, $\qe x\,A(x)$, respectively) itself.  

  The set of truth-values of subformulas of $A$ under a given
  interpretation~$\I$ is denoted by
  \[ \subval(\I,A) = \{ \I(B) \suchthat B 
  \text{ subformula of $A$ w.r.t.\ $U^\I$}\} \cup \{0,1\}
   \] 
   If $\Gamma$ is a set of formulas, then $\subval(\I, \Gamma) = \bigcup
   \{\subval(\I, A) : A \in \Gamma\}$.
\end{definition}

\begin{lemma}\label{lem:induced-interpretation}
  Let $\I$ be a $V$-interpretation, and let $h\colon \subval(\I, \Gamma) \to
  W$ be a mapping satisfying the following properties:
\begin{enumerate}
\item $h(0) = 0$, $h(1) = 1$;
\item $h$ is strictly monotonic, i.e., if $a < b$, then $h(a) < h(b)$;
\item for every $X \subseteq \subval(\I, \Gamma)$, $h(\inf X) = \inf h(X)$ and
$h(\sup X) = \sup h(X)$ (provided $\inf X$, $\sup X \in \subval(\I, \Gamma)$).
\end{enumerate}
Then the $W$-interpretation $\I_h$ 
with universe $U^\I$, $f^{\I_h} = f^\I$, and for atomic $B \in\LL^\I$,
\[
    \I_h(B) = \begin{cases} h(\I(B)) & \text{ if $\I(B) \in \dom h$} \\ 
                            1 &        \text{ otherwise}
               \end{cases}
\]
satisfies $\I_h(A) = h(\I(A))$ for all $A \in \Gamma$.
\end{lemma}

\begin{proof}
  By induction on the complexity of~$A$.  If $A \equiv \bot$, the claim
  follows from (1).  If $A$ is atomic, it follows from the definition of
  $\I_h$.  For the propositional connectives the claim follows from the strict
  monotonicity of~$h$ (2).  For the quantifiers, it follows from property~(3).
\end{proof}

\begin{remark}
  Note that the construction of $\I_h$ and the proof of
  Lemma~\ref{lem:induced-interpretation} also goes through without the
  condition $h(0) = 0$, provided that the formulas in $\Gamma$ do not
  contain~$\bot$, and goes through without the requirement that existing inf's
  be preserved ($h(\inf X) = \inf h(X)$ if $\inf X \in \subval(\I, \Gamma)$)
  provided they do not contain~$\forall$.
\end{remark}

\begin{definition}
  A \emph{\g-embedding} $h\colon V \to W$ is a strictly monotonic, continuous
  mapping between Gödel sets which preserves
  $0$ and $1$.
\end{definition}

\begin{lemma}\label{lem:g-embed}
  Suppose $h\colon V \to W$ is a \g-embedding. (a) If $\I$ is a
  $V$-interpretation, and $\I_h$ is the interpretation induced by $\I$ and
  $h$, then $\I_h(A) = h(\I(A))$. (b) If $\Gamma \entails_W A$ then $\Gamma
  \entails_V A$ (and hence $\gdl{W} \subseteq \gdl{V}$). (c) If $h$ is
  bijective, then $\Gamma \entails_W A$ iff $\Gamma \entails_V A$ (and hence,
  $\gdl{V} = \gdl{W}$).
\end{lemma}

\begin{proof}
  (a) $h$ satisfies the conditions of Lemma~\ref{lem:induced-interpretation},
  for $\Gamma$ the set of all formulas. (b) If $\Gamma \nentails_V A$, then
  for some $\I$, $\I(B) = 1$ for all $B \in \Gamma$ and $\I(A) < 1$.  By
  Lemma~\ref{lem:induced-interpretation}, $\I_h(B) = 1$ for all $B \in \Gamma$
  and $\I_h(A) < 1$ (by strict monotonicity of $h$). Thus $\Gamma \nentails_W
  A$. (c) If $h$ is bijective then $h^{-1}$ is also a \g-embedding.
\end{proof}

\begin{definition}[Submodel, elementary submodel]
  Let $\I_1$, $\I_2$ be interpretations.  We write $\I_1 \subseteq \I_2$ 
  ($\I_2$ extends $\I_1$) iff
  $U^{\I_1} \subseteq U^{\I_2}$, and for all $k$, all $k$-ary 
  predicate symbols $P$ in $\LL$, and all $k$-ary function symbols $f$
  in $\LL$ we have
  \[ P^{ \I_1 } =  P^{ \I_2 } \upharpoonright (U^{\I_1})^k  \quad
     f^{ \I_1 } =  f^{ \I_2 } \upharpoonright (U^{\I_1})^k  \quad
     \]
  or in other words, if $\I_1$ and $\I_2$ agree on closed atomic
  formulas.  

  We write $\I_1 \prec \I_2$ if $\I_1 \subseteq \I_2$ and $\I_1(A) = \I_2(A)$ 
  for all $\LL^{U^{\I_1}}$-formulas $A$. 
\end{definition}

\begin{proposition}[Downward Löwenheim-Skolem]\label{fact:loewenheim}
  For any interpretation $\I$ with $U^\I$ infinite, there is an
  interpretation $\I' \prec \I$  with a countable universe~$U^{\I'} $. 
\end{proposition}

\begin{proof}[Proof sketch]
  The proof is an easy generalization of the construction for the classical
  case.  We construct a sequence of countable subsets $U_1 \subseteq U_2
  \subseteq \cdots$ of $U^\I$: $U_1$ simply contains $t^\I$ for all closed
  terms of the original language.  $U_{i+1}$ is constructed from $U_i$ by
  adding, for each of the (countably many) formulas of the form $\exists x\,
  A(x)$ and $\forall x\, A(x)$ in the language $\LL^{U_i}$, a countable
  sequence $a_j$ of elements of $U^\I$ so that $(\I(A(a_j)))_j \to \I(\exists
  x\, A(x))$ or $\to \I(\forall x\, A(x))$, respectively.  $U^{\I'} =
  \bigcup_i U_i$.
\end{proof}

\begin{lemma}\label{lemma:interpretation-cut-off}
  Let $\I$ be a interpretation into $V$, $w\in[0,1]$, and let $\I_w$
  be defined by
  \[ \I_w(B) = \begin{cases} \I(B) &\text{ if $\I(B)<w$}\\
    1 &\text{ otherwise}
                \end{cases}
  \]
  for atomic formulas $B$ in $\LL^U$.  Then $\I_w$ is an interpretation into
  $V$.  If $w \notin \subval(\I, A)$, then $\I_w(A) = \I(A)$ if $\I(A) < w$,
  and $\I_w(A) = 1$ otherwise.
\end{lemma}

\begin{proof}
  Let $h_w(a)=a$ if $a<w$ and $=1$ otherwise.  By induction on the complexity
  of formulas $B$ it is easily shown that $\I'(B)=h_w(\I(B))$ for all
  subformulas $B$ of $A$ w.r.t.\ $U^\I$.
\end{proof}

\begin{proposition}\label{prop:consequences}
$\Gamma \entails A$ iff $\Gamma \Vdash A$
\end{proposition}

\begin{proof}
  Only if: obvious. If: Suppose that $\Gamma \nentails A$, i.e., there is a
  $V$-interpretation $\I$ so that $\inf \{\I(B) : B \in \Gamma\} > \I(A)$.  By
  Proposition~\ref{fact:loewenheim}, we may assume that $U^\I$ is countable.
  Hence, there is some $w$ with $\I(A) < w < \inf\{\I(B) : B \in \Gamma\}$ and
  $w \notin \subval(\I, \Gamma \cup \{A\})$.  Let $\I_{w}$ be as in
  Lemma~\ref{lemma:interpretation-cut-off}.  Then $\I_w(B) = 1$ for all $B \in
  \Gamma$ and $\I_w(A) < 1$.
\end{proof}

The coincidence of the two consequence relations is a unique feature of Gödel
logics.  Proposition~\ref{prop:consequences} does not hold in \L ukasiewicz
logic, for instance. There, $A, A \limp_\textrm{\L} B \Vdash B$ but $A, A
\limp_\textrm{\L} B \nentails B$.  In what follows, we will use $\entails$
when semantic consequence is at issue; the preceding propositions shows that
the results we obtain for $\entails$ hold for $\Vdash$ as well.

\begin{lemma}[Semantic deduction theorem]\label{thm:semdeduc}
  \[ \Gamma, A \entails B \quad\text{iff}\qquad \Gamma\entails A\limp B.\]
\end{lemma}

\begin{proof}
  Immediate consequence of the definition of $\entails$ and the
  semantics for $\limp$.
\end{proof}

We want to conclude this part with two interesting observations:

\paragraph{Relation to residuated algebras}
If one considers the truth value set as a Heyting algebra with $a\wedge
b = \min(a,b)$,  $a\vee b = \max(a,b)$, and
\[ a \to b = \left\{
  \begin{array}{ll}
    1 & \mbox{if $a\le b$ }\\
    b & \mbox{otherwise}\\
  \end{array}\right.
\]
then $\to$ and $\wedge$ are residuated, i.e.,
\[ 
  (a \to b)  = \sup \{\, x: \, (x\wedge a)\, \le\, b\, \}. 
\]

\paragraph{The Gödel conditional}
A large class of many-valued logics can be developed from the theory of
$t$-norms \cite{hajek}. The class of $t$-norm based logics includes not only
(standard) Gödel logic, but also {\L}ukasiewicz- and product logic. In these
logics, the conditional is defined as the residuum of the respective $t$-norm,
and the logics differ only in the definition of their $t$-norm and the
respective residuum, i.e., the conditional. The truth function for the Gödel
conditional is of particular interest as it can be `deduced' from simple
properties of the evaluation and the entailment relation, a fact which was
first observed by G.~Takeuti.
\begin{lemma}\label{lemma:takeuti}
  Suppose we have a standard language containing a `conditional'
  $\twoheadrightarrow$ interpreted by a truth-function into $[0,1]$. Suppose
  further that 
  \begin{enumerate}
  \item a conditional evaluates to~$1$ if the truth value of the antecedent
    is less or equal to the truth value of the consequent, i.e.,
    if $\I(A) \le \I(B)$, then $\I(A \twoheadrightarrow B) = 1$;
  \item $\entails$ is defined as above, i.e., if $\Gamma \entails B$, then
    $\min \{ \I(A) \suchthat A\in\Gamma \} \le \I(B)$;
  \item the deduction theorem holds, i.e., $\Gamma \cup \{ A \} \entails B
    \dequiv \Gamma \entails A \twoheadrightarrow B$.
  \end{enumerate}
  Then $\twoheadrightarrow$ is the Gödel conditional.
\end{lemma}

\begin{proof}
  From (1), we have that $\I(A \twoheadrightarrow B) = 1$ if $\I(A) \le
  \I(B)$.  Since $\entails$ is reflexive, $B \entails B$. Since it is
  monotonic, $B, A \entails B$.  By the deduction theorem, $B \entails A
  \twoheadrightarrow B$.  By (2),
  \[ 
  \I(B) \le \I(A \twoheadrightarrow B).
  \]
  From $A \twoheadrightarrow B \entails A \twoheadrightarrow B$ and the
  deduction theorem, we get $A \twoheadrightarrow B, A \entails B$. By (2),
  \[ 
  \min \{\I(A \twoheadrightarrow B), \I(A)\} \le \I(B).
  \]
  Thus, if $\I(A) > \I(B)$, $\I(A \twoheadrightarrow B) \le \I(B)$.
\end{proof}

Note that all usual conditionals (Gödel, \L ukasiewicz, product conditionals)
satisfy condition (1).  So, in some sense, the Gödel conditional is the only
many-valued conditional which validates both directions of the deduction
theorem for $\entails$.  For instance, for the \L ukasiewicz conditional
$\limp_\textrm{\L}$ the right-to-left direction fails: $A \limp_\textrm{\L} B
\entails A \limp_\textrm{\L} B$, but $A \limp_\textrm{\L} B, A
\nentails B$. (With respect to $\Vdash$, the left-to-right direction of the
deduction theorem fails for $\limp_\textrm{\L}$.)


\subsection{Axioms and deduction systems}

In this section we introduce certain axioms and deduction systems for
Gödel logics, and we will show completeness of these deduction systems
subsequently. We will use a Hilbert style proof system:
\begin{definition}
  A formula~$A$ is derivable from formulas $\Gamma$ in a system~$\cA$
  consisting of the axioms and the rules iff there are formulas $A_0$, \dots,
  $A_n = A$ such that for each $0 \le i \le n$ either $A_i \in \Gamma$, or
  $A_i$ is an instance of an axiom in $\cA$, or there are indices $j_1$,
  \dots, $j_l < i$ and a rule in $\cA$ such that $A_{j_1}$, \dots, $A_{j_l}$
  are the premises and $A_i$ is the conclusion of the rule. In this case we
  write $\Gamma \proves_{\cA} A$.
\end{definition}

We will denote by $\IL$ the following complete axiom system for intuitionistic
logic (taken from \cite{troelstra-handbook}). Rules are written as $A_1,
\ldots, A_n \proves A$.  \def\ax(#1){\hbox{({#1})}} 
  \fxnote{Besseres Axiomensystem mit weniger Regeln finden!}
\begin{align*}
  \ax(I1) &\quad A, A\limp B\proves B &
  \ax(I2) &\quad A\limp B,B\limp C\proves A\limp C\\
  \ax(I3) &\quad A\lor A\limp A, A\limp A\land A  &
  \ax(I4) &\quad A\limp A\lor B, A\land B\limp A\\
  \ax(I5) &\quad A\lor B\limp B\lor A, A\land B\limp B\land A &
  \ax(I6) &\quad A\limp B\proves C\lor A\limp C\lor B\\
  \ax(I7) &\quad A\land B\limp C \proves A\limp(B\limp C) &
  \ax(I8) &\quad A\limp (B\limp C) \proves A\land B\limp C\\
  \ax(I9) &\quad \bot\limp A \\
  \ax(I10) &\quad B^{(x)}\limp A(x) \proves B^{(x)}\limp\qa xA(x) &
  \ax(I11) &\quad \qa xA(x) \limp A(t)\\
  \ax(I12) &\quad A(t)\limp \qe xA(x) &
  \ax(I13) &\quad A(x)\limp B^{(x)} \proves \qe xA(x)\limp B^{(x)}
\intertext{(where $B^{(x)}$ means that $x$ is not free in $B$).}
\end{align*}
The following axioms will play an important r\^ole ($\axqs$ stands for
`quantifier shift', $\axlin$ for `linearity', $\axiso$ for `isolation
axiom of $0$', and $\axfinn$ for `finite with $n$ elements'):
\begin{align*}
  \axqs  &\quad \qa x(C^{(x)}\lor A(x)) \limp (C^{(x)}\lor\qa xA(x)) \\
  \axlin &\quad (A\limp B)\lor (B\limp A)\\
  \axiso &\quad \qa x\lnot\lnot A(x) \limp  \lnot\lnot\qa xA(x)\\
  \axfinn &\quad (\top\limp A_1) \lor (A_1\limp A_2)\lor \ldots \lor
    (A_{n-2}\limp A_{n-1})\lor (A_{n-1}\limp\bot)
\end{align*}

\begin{notation}
  $\H$ denotes the axiom system $\IL + \axqs + \axlin$.

  $\H_n$ for $n \ge 2$ denotes the axiom system $\H + \axfinn$.

  $\Ho$ denotes the axiom system $\H + \axiso$.
\end{notation}

\begin{theorem}[Soundness]\label{thm:soundness}
  Suppose $\Gamma$ contains only closed formulas, and all axioms of $\cA$ are
  valid in $\gdl{V}$. Then, if $\Gamma \proves_\cA A$ then $\Gamma \entails_V
  A$.
\end{theorem}

\begin{proof}
  By induction on the complexity of proofs.  By assumption, all axioms
  of~$\cA$ are valid in $\gdl{V}$, hence $\Gamma \entails_V A_i$ if $A_i$ is
  an axiom.  If $A_i \in \Gamma$, then obviously $\Gamma \entails_V A_i$. 
  It remains to show that the rules of inference preserve consequence.  We
  show this for modus ponens (I1) and existential generalization (I13), the
  other cases are analogous.

  Suppose $\Gamma \entails_V A$ and $\Gamma \entails_V A \limp B$ and consider
  a $V$-interpretation $\I$.  Let $v = \inf \{\I(C) : C \in \Gamma\}$.  If
  $\I(A) \le \I(B)$, then we have $v \le \I(B)$ because $v \le \I(A)$.  If
  $\I(A) > \I(B)$, then $v \le \I(B)$ because $\I(B) = \I(A \limp B)$.
  
  Suppose $\Gamma \entails_V A(x) \limp B$ and $x$ does not occur free in $B$.
  Let $\I$ be a $V$-interpretation, and let $w = \sup \{\I(A(u)) : u \in
  U^\I\}$, and let $\I_u$ be the interpretation resulting from $\I$ by
  assigning $u$ to $x$.  Since the formulas in $\Gamma$ are all closed and $B$
  does not contain $x$ free, $\I_u(C) = \I(C)$ for all $C \in \Gamma \cup
  \{B\}$ and $u \in U^\I$. Now suppose $w > \I(\exists x\, A(x) \limp B)$. In
  this case, $\I(\exists x\, A(x)) > \I(B)$.  But then, for some $u \in U^\I$,
  $\I_u(A(x)) > \I(B)$ and we'd have $w > \I_u(A(x) \limp B)$, contradicting
  $\Gamma \entails_V A(x) \limp B$.  The case for (I10) is analogous.
\end{proof}

Note that the restriction to closed formulas in $\Gamma$ is essential: $A(x)
\proves_\H \forall x\, A(x)$ but obviously $A(x) \nentails_V \forall x\,
A(x)$.

\subsection{Relationships between Gödel logics}\label{sec:intro:relgoedel}

The relationships between finite and infinite valued \emph{propositional}
Gödel logics are well understood.  Any choice of an infinite set of
truth-values results in the same propositional Gödel logic, viz., Dummett's
\LC.  \LC~was defined using the set of truth-values~$V_\dn$ (see below).
Furthermore, we know that \LC{} is the intersection of all finite-valued
propositional Gödel logics, and that it is axiomatized by intuitionistic
propositional logic~\IPL{} plus the schema $(A \impl B) \lor (B \impl A)$.
\IPL{} is contained in all Gödel logics.

In the first-order case, the relationships are somewhat more
interesting. First of all, let us note the following fact
corresponding to the end of the previous paragraph:

\begin{proposition}
Intuitionistic predicate logic~\IL{} is contained in all first-order
Gödel logics.
\end{proposition}

\begin{proof}
The axioms and rules of \IL{} are sound for the Gödel truth
functions.
\end{proof}

As a consequence of this proposition, we will be able to use any
intuitionistically sound rule and intuitionistically true formula
when working in any of the Gödel logics.

We can consider special truth value sets which will act as
prototypes for other logics. This is due to the fact that the logic is
defined extensionally as the set of formulas valid in this truth value
set, so the Gödel logics on different truth value sets may coincide.
\begin{align*}
V_\R & =  [0,1] \\
V_\vdn & =  \{1/k : k \ge 1\} \cup \{0\} \\
V_\vup & =  \{1 - 1/k : k \ge 1\} \cup \{1\} \\
V_m & = \{1 - 1/k : 1 \le k \le m-1\} \cup \{1\}
\end{align*}
The corresponding Gödel logics are $\gdl \R$, $\gdl \vdn$, $\gdl \vup$,
$\gdl m$. $\gdl \R$ is the \emph{standard} Gödel logic.

The logic $\gdl \vdn$ also turns out to be closely related 
to some temporal logics \cite{BaazLeitschZach:96,BaazLeitZach96TCS}.
$\gdl \vup$ is the intersection of all finite-valued first-order Gödel
logics as shown in Theorem~\ref{herbrand}.

\fxnote{Sollen die Ergebnisse dieses Abschnitts alle für $\entails_V$ statt für
  $\gdl V$ formuliert werden?}

\begin{proposition}
$\gdl \R = \bigcap_V \gdl V$, where $V$ ranges over all Gödel sets.
\end{proposition}

\begin{proof}
If $\gdl V \entails A$ for every $V$, then also for $V = [0, 1]$.  Conversely,
if there is some Gödel set $V$ and a $V$-interpretation $\I$ with $\I \nmodels
A$, then $\I$ is also a $[0,1]$-interpretation and hence $\gdl \R \nentails A$.
\end{proof}

\begin{proposition}\label{basic-cont}
The following strict containment relationships hold:
\begin{enumerate}
\item $\gdl m \supsetneq \gdl {m+1}$,
\item $\gdl m \supsetneq \gdl \vup \supsetneq \gdl \R$,
\item $\gdl m \supsetneq \gdl \vdn \supsetneq \gdl \R$.
\end{enumerate}
\end{proposition}

\begin{proof}
The only non-trivial part is proving that the containments are strict.
For this note that
\[
(A_1 \impl A_2) \lor \ldots \lor (A_m \impl A_{m+1})
\]
is valid in $\gdl m$ but not in $\gdl {m+1}$. Furthermore, let
\begin{align*}
C_\vup & =  \exists x(A(x) \impl \forall y\, A(y)) \mathrm{\ and} \\
C_\vdn & =  \exists x(\exists y\, A(y) \impl A(x)).
\end{align*}
$C_\dn$ is valid in all $\gdl m$ and in $\gdl \vup$ and $\gdl \dn$; $C_\vup$ is
valid in all $\gdl m$ and in $\gdl \vup$, but not in $\gdl \dn$; neither is
valid in $\gdl \R$ (\cite{BaazLeitschZach:96}, Corollary~2.9).
\end{proof}

The formulas $C_\vup$ and $C_\dn$ are of some importance in the study
of first-order infinite-valued Gödel logics.  $C_\vup$ expresses the
fact that every infimum in the set of truth values is a minimum, and
$C_\dn$ states that every supremum (except possibly $1$) is a maximum.
The intuitionistically admissible quantifier
shifting rules are given by the following implications and equivalences:
\[
\begin{array}{rcl}
(\forall x\, A(x) \land B) & \equiv & \forall x(A(x) \land B) \\
(\exists x\, A(x) \land B) & \equiv & \exists x(A(x) \land B) \\
(\forall x\, A(x) \lor B) & \impl & \forall x(A(x) \lor B) \\
(\exists x\, A(x) \lor B) & \equiv & \exists x(A(x) \lor B) \\
(B \impl \forall x\, A(x)) & \equiv & \forall x(B \impl A(x)) \\
(B \impl \exists x\, A(x)) & \rimpl & \exists x(B \impl A(x)) \\
(\forall x\, A(x) \impl B) & \rimpl & \exists x(A(x) \impl B) \\
(\exists x\, A(x) \impl B) & \equiv & \forall x(A(x) \impl B)
\end{array}
\]
The remaining three are:
\[
\begin{array}{rcl}
(\forall x\, A(x) \lor B) & \rimpl & \forall x(A(x) \lor B) \\
(B \impl \exists x\, A(x)) & \impl & \exists x(B \impl A(x)) \\
(\forall x\, A(x) \impl B) & \impl & \exists x(A(x) \impl B) 
\end{array}
\eqno{\begin{array}{r}(S_1)\\(S_2)\\(S_3)\end{array}}
\]
Of these, $S_1$ is valid in any Gödel logic.  $S_2$ and
$S_3$ imply and are implied by $C_\dn$ and $C_\vup$, respectively
(take $\exists y\, A(y)$ 
and $\forall y\, A(y)$, respectively, for $B$).  $S_2$ and $S_3$ are,
respectively, both valid in $\gdl \vup$, invalid and valid in $\gdl \dn$, and
both invalid in $\gdl \R$. Thus we obtain
\begin{corollary}
  $\gdl \vup$ is the only Gödel logic where every formula is equivalent to
  a prenex formula with the same propositional matrix.
\end{corollary}

We now also know that $\gdl \vup \neq \gdl \dn$.  In fact, we have
$\gdl \dn \subsetneq \gdl \vup$; this follows from the following theorem.

\begin{theorem}\label{herbrand}
\[\gdl \vup = \bigcap\limits_{m \ge 2} \gdl m\]
\end{theorem}

\begin{proof}
  By Proposition~\ref{basic-cont}, 
  $\gdl \vup \subseteq \bigcap_{m \ge 2}\gdl m$.  We now prove the reverse
  inclusion. Assume that there is an interpretation $\I$ such that
  $\I\nmodels A$, we want to give an interpretation $\I'$ such that
  $\I'\nmodels A$ and $\I'$ is a $\gdl m$ interpretation for
  some~$m$.
  
  Suppose there is an interpretation $\I$ such that $\I\nmodels A$, let $\I(A)
  = 1 - 1/k$. Let $w$ be somewhere between $1-1/k$ and $1-1/(k+1)$. Then the
  interpretation $\I_w$ given in Lemma~\ref{lemma:interpretation-cut-off} also
  is a counterexample for $A$.  Since there are only finitely many truth
  values below $w$ in $V_\vup$, $\I_w$ is a $\gdl {k+1}$ interpretation with
  $\I_w \nmodels A$. This completes the proof of the theorem.
\end{proof}

\begin{corollary}
$\gdl m \supsetneq \bigcap_m  \gdl m = 
   \gdl \vup \supsetneq \gdl \dn \supsetneq \gdl \R$
\end{corollary}

As we will see later, the axioms $\axfinn$ axiomatize exactly the
finite-valued Gödel logics. In these logics the quantifier shift axiom $\axqs$
is not necessary. Furthermore, all quantifier shift rules are valid in the
finite valued logics. Since $\gdl \vup$ is the intersection of all the finite
ones, all quantifier shift rules are valid in~$\gdl \vup$.  Moreover, any
infinite-valued Gödel logic other than $\gdl \vup$ is defined by some~$V$ which
either contains an infimum which is not a minimum, or a supremum (other
than~$1$) which is not a maximum.  Hence, in $V$ either $C_\vup$ or $C_\vdn$
will be invalid, and therewith either $S_3$ or $S_2$. We have:

\begin{corollary}
  $\gdl \vup$ is the only Gödel logic with infinite truth value set which
  admits all quantifier shift rules.
\end{corollary}

\section{Topology and Order}\label{sec:order}

%
%
\subsection{Perfect sets}

All the following notations, lemmas, theorems are carried out within the
framework of Polish spaces, which are separable, completely metrizable
topological spaces. For our discussion it is only necessary to know
that~$\bbR$ and all its closed subsets are Polish spaces (hence, every Gödel
set is a Polish space). For a detailed exposition see
\cite{Moschovakis:1980,kechris}.

\begin{definition}[limit point, perfect space, perfect set]
        \hskip0pt plus3pt minus 3pt
  A \emph{limit point} of a topological space is a point that is not
  isolated, i.e.\ for every open neighborhood~$U$ of~$x$ there is a
  point~$y\in U$ with~$y\neq x$. A space is \emph{perfect} if all its
  points are limit points. A set $P\subseteq \bbR$ is \emph{perfect}
  if it is closed and together with the topology induced from $\bbR$
  is a perfect space.
\end{definition}

It is obvious that all (non-trivial) closed intervals are perfect sets, also
all countable unions of (non-trivial) intervals. But all these sets generated
from closed intervals have the property that they are `everywhere dense',
i.e., contained in the closure of their inner component. There is another very
famous set which is perfect but is nowhere dense, the Cantor set:

\begin{example}[Cantor Set]
  The set of all numbers in the unit interval which can be expressed in
  triadic notation only by digits~$0$ and~$2$ is called \emph{Cantor
    set}~$\bbD$.
\end{example}

A more intuitive way to obtain this set is to start with the unit
interval, take out the open middle third and restart this process with
the lower and the upper third. Repeating this you get exactly the
Cantor set because the middle third always contains the numbers which
contain the digit $1$ in their triadic notation.

This set has a lot of interesting properties, the most important one
is that it is a perfect set:

\begin{proposition}
  The Cantor set is perfect.
\end{proposition}

It is possible to embed the Cauchy space into any perfect space, yielding
the following proposition:

\begin{proposition}\label{lem:uncountable}
  If~$X$ is a nonempty perfect Polish space, then the cardinality
  of~$X$ is~$2^{\alepho}$ and therefore, all nonempty perfect subsets,
  too, have cardinality of the continuum.
\end{proposition}

It is possible to obtain the following characterization of perfect
sets (see \cite{winkler:howmuch}):

\begin{proposition}[Characterization of perfect sets in $\bbR$]\label{prop:winkler}
  For any perfect subset of $\bbR$ there is a unique partition of the
  real line into countably many intervals such that the intersections
  of the perfect set with these intervals are either empty, the full
  interval or isomorphic to the Cantor set.
\end{proposition}

So we see that intervals and Cantor sets are prototypical for perfect
sets and the basic building blocks of more complex perfect sets.

Every Polish space can be partitioned into a perfect kernel and a
countable rest. This is the well known Cantor-Bendixon Theorem:

%

\begin{theorem}[Cantor-Bendixon]\label{thm:cantorbendixon}
  Let~$X$ be a Polish space. Then~$X$ can be uniquely written as $X =
  P \cup C$, with~$P$ a perfect subset of~$X$ and~$C$ countable and
  open. The subset~$P$ is called the \emph{perfect kernel} of~$X$
  (denoted with $V^\infty$).
\end{theorem}

As a corollary we obtain that any uncountable Polish space contains a
perfect set, and therefore, has cardinality~$2^{\alepho}$.

\subsection{Relation to Gödel logics}


The following lemma was originally proved in \cite{Prei03PHD}, where it was
used to extend the proof of recursive axiomatizability of `standard' Gödel
logics (those with $V=[0,1]$) to Gödel logics with a truth value set
containing a perfect set in the general case. The following more
simple proof is inspired by \cite{bgp}:





\begin{lemma}\label{lemma:countable-into-cantor}
  Suppose that $M \subseteq [0,1]$ is countable and $P \subseteq [0,1]$ is
  perfect. Then there is a strictly monotone continuous map 
  $h\colon M \to  P$ (i.e., infima and suprema already existing in $M$
  are preserved). Furthermore, if $\inf M \in M$, then one can choose
  $h$ such that $h(\inf M) = \inf P$.
\end{lemma}

\begin{proof}
  Let $\sigma$ be the mapping which scales and shifts $M$ into $[0,1]$,
  i.e.\ the mapping $x\to (x-\inf M)/(\sup M-\inf M)$ (assuming that
  $M$ contains more than one point).
  Let $w$ be an injective monotone map from $\sigma(M)$ into $2^\omega$,
  i.e. $w(m)$ is a fixed binary representation of $m$. For dyadic rational
  numbers (i.e.\ those with different binary representations) we fix one
  possible.
  
  Let $i$ be the natural bijection from $2^\omega$ (the set of infinite
  $\{0,1\}$-sequences, ordered lexicographically) onto $\bbD$, the Cantor set.
  $i$ is an order preserving homeomorphism.  Since $P$ is perfect, we can find
  a continuous strictly monotone map $c$ from the Cantor set $\bbD \subseteq
  [0,1]$ into $P$, and $c$ can be chosen so that $c(0) = \inf P$.
  
  Now $h = c \circ i \circ w \circ \sigma$ is also a strictly monotone
  map from $M$ into~$P$, and $h(\inf M) = \inf P$,  if $\inf M \in M$.
  Since $c$ is continuous, existing infima and suprema are preserved.
\end{proof}

\begin{corollary}\label{bendixon}
  A Gödel set~$V$ is uncountable iff it contains a non-trivial dense
  linear subordering.
\end{corollary}

\begin{proof}
  If: Every countable non-trivial dense linear order has order type
  $\boeta$, $\one+\boeta$, $\boeta+\one$, or $\one+\boeta+\one$
  \cite[Corollary~2.9]{rosenstein}, where $\boeta$ is the order type
  of $\bbQ$.  The completion of any ordering of order type~$\boeta$
  has order type~$\lambda$, the order type of $\bbR$
  \cite[Theorem~2.30]{rosenstein}, thus the truth value set must be
  uncountable.
  
  Only if: By Theorem~\ref{thm:cantorbendixon}, $V^\infty$ is non-empty. Take
  $M = \bbQ \cap [0, 1]$ and $P = V^\infty$ in
  Lemma~\ref{lemma:countable-into-cantor}.  The image of $M$ under~$h$ is a
  non-trivial dense linear subordering in~$V$.
\end{proof}

\begin{theorem}\label{thm:gs-id}
  Suppose $V$ is a truth value set with non-empty perfect kernel $P$, and
  let $W = V \cup [\inf P,1]$. Then 
  $\mathnormal{\entails_V} = \mathnormal{\entails_W}$, i.e.\
  $\Gamma \entails_V A$ iff $\Gamma \entails_W A$. Thus also the logics
  induced by $V$ and $W$ are the same, i.e., $\gdl{V} = \gdl{W}$.
\end{theorem}

\begin{proof}
  As $V \subseteq W$ we have 
  $\mathnormal{\entails_W} \subseteq \mathnormal{\entails_V}$ 
  (cf.\ the Remark preceding Definition~\ref{def:goedellogics}).
  Now assume that $\I$ is a $W$-interpretation which shows that
  $\Gamma\entails_W A$ does \emph{not} hold, i.e.,
  $\inf\{\I(B)\suchthat B\in\Gamma\} > \I(A)$. By
  Proposition~\ref{fact:loewenheim}, we may assume that $U^\I$ is countable.
  The set $\subval(\I, \Gamma\cup A)$ has cardinality at most
  $\aleph_0$, thus there is a $b \in [0, 1]$ such that 
  $b \notin \subval(\I, \Gamma\cup A)$ and $\I(A) < b < 1$.  By 
  Lemma~\ref{lemma:interpretation-cut-off}, $\I_b(A) < b < 1$.  Now consider
  $M = \subval(\I_b, \Gamma\cup A)$: these are all the truth values
  from $W = V \cup [\inf P, 1]$ required to compute $\I_b(A)$ and
  $\I_b(B)$ for all $B\in\Gamma$.  We have to find some way to map
  them to $V$ so that the induced interpretation is a counterexample
  to~$\Gamma\entails_V A$.  

  Let $M_0 = M \cap [0, \inf P)$ and 
  $M_1 = (M \cap [\inf P, b])\cup\{\inf P\}$.
  By Lemma~\ref{lemma:countable-into-cantor} there is a strictly monotone
  continuous (i.e. preserving all existing infima and suprema) map~$h$
  from $M_1$ into $P$. Furthermore, we can choose~$h$ such that
  $h(\inf M_1) = \inf P$.

  We define a function $g$ from $\subval(\I_b,\Gamma\cup A)$ to $V$ as
  follows:
  \[
     g(x) = \begin{cases}
               x          & 0\le x \le\inf P\\
               h(x)       & \inf P \le x \le b\\
               1          & x = 1
            \end{cases}
  \]
  Note that there is no $x\in\subval(\I_b,\Gamma\cup A)$ with $b<x<1$.
  This function has the following properties: $g(0)=0$, $g(1)=1$, $g$
  is strictly monotonic and preserves existing infima and
  suprema. Using Lemma~\ref{lem:induced-interpretation} we obtain that
  $\I_g$ is a $V$-interpretation with $\I_g(C) = g(\I_b(C))$ for all
  $C\in\Gamma\cup A$, thus also 
  $\inf\{\I_g(B)\suchthat B\in\Gamma\} > \I_g(A)$.
\end{proof}

\section{Countable Gödel sets}\label{sec:countable}

In this section we show that the first-order Gödel logics where the set of
truth values does not contain a dense subset are not axiomatizable. We
establish this result by reducing the classical validity of a formula in all
finite models to the validity of a formula in Gödel logic (the set of these
formulas is not r.e.\ by Trakhtenbrot's Theorem).

\begin{definition}
  A formula is called \emph{crisp} if all occurrences of atomic
  formulas are either negated or double-negated.
\end{definition}

\begin{lemma}\label{lm:crisp} If $A$ and $B$ are crisp and classically
  equivalent, then also $\GR \models A \leftrightarrow B$. Specifically,
  if $A(x)$ and $B$ are crisp, then
  \begin{align*}
    & \entails \qa xA(x)\limp B \leftrightarrow
    \qe x(A(x)\limp B) \quad\text{and} \\
    & \entails B\limp \qe xA(x) \leftrightarrow
    \qe x(B\limp A(x)).
  \end{align*}
\end{lemma}

\begin{proof}
  Given an interpretation~$\I$, define $\I'(C) = 1$ if $\I(C) > 0$ and $= 0$
  if $\I(C) = 0$ for atomic~$C$.  It is easily seen that if $A$, $B$ are
  crisp, then $\I(A) = \I'(A)$ and $\I(B) = \I'(B)$.  But $\I'$ is a classical
  interpretation, so by assumption $\I'(A) = \I'(B)$.
\end{proof}

\begin{theorem}\label{thm:count-nonax}
  If $V$ is countably infinite, then $\gdl V$ is not recursively
  enumerable. 
\end{theorem}

\begin{proof}
By Theorem~\ref{bendixon}, $V$ is countably infinite iff it is
infinite and does not contain a non-trivial densely ordered subset.
We show that for every sentence $A$ there is a sentence $A^g$
s.t. $A^g$ is valid in $\gdl V$ iff $A$ is true in every finite
(classical) first-order structure.

We define $A^g$ as follows: Let $P$ be a unary and $L$ be a binary
predicate symbol not occurring in~$A$ and let $Q_1$, \dots, $Q_n$ be
all the predicate symbols in~$A$.  We use the abbreviations $x \in y
\equiv \neg\neg L(x, y)$ and $x \prec y \equiv (P(y) \impl P(x)) \impl
P(y)$.  Note that for any interpretation~\I, $\I(x \in y)$ is either
$0$ or $1$, and as long as $\I(P(x)) < 1$ for all $x$ (in particular,
if $\I(\exists z\, P(z)) < 1$), we have $\I(x \prec y) = 1$ iff
$\I(P(x)) < \I(P(y))$.  Let $A^g \equiv$
\begin{equation}\label{fm:countable}
  \left\{
    \begin{array}{l}
      S \land {} c_1 \in 0 \land {} c_2 \in 0 \land 
      c_2 \prec c_1 \land {}\\
      \quad \forall i\bigl[\forall x, y \forall j \forall k \exists z\, D
      \lor \forall x\neg(x \in s(i))\bigr]
    \end{array}\right\}
  \impl (A' \lor \exists u\, P(u))
\end{equation}
where $S$ is the conjunction of the standard axioms for $0$, successor
and $\le$, with double negations in front of atomic formulas,
\[ D\equiv \begin{array}{l}
(j \le i \land x \in j \land k\le i \land y \in k \land x \prec y) \impl {}\\
\qquad \impl (z \in s(i) \land x \prec z \land z \prec y)
\end{array}\]
and
$A'$ is $A$ where every atomic formula is replaced by its double
negation, and all quantifiers are relativized to the predicate
$R(i)\equiv \exists x(x \in i)$.

Intuitively, $L$ is a predicate that divides a subset of the domain
into levels, and $x \in i$ means that $x$ is an element of level~$i$.
If the antecendent is true, then the true standard axioms~$S$ force
the domain to be a model of PA, which could be either a standard model
(isomorphic to $\bbN$) or a non-standard model ($\bbN$ followed by copies
of $\bbZ$).
$P$ orders the elements of the domain which fall into one of the
levels in a subordering of the truth values.  

The idea is that for any
two elements in a level $\le i$ there is an element in a not-empty
level $j\ge i$
which lies strictly between those two elements in the ordering given
by~$\prec$.  If this condition cannot be satisfied, the levels above
$i$ are empty.  Clearly, this condition can be satisfied in an
interpretation~\I{} only for finitely many levels if $V$ does not
contain a dense subset, since if more than finitely many levels are
non-empty, then $\bigcup_{i} \{\I(P(d)) : \I \models d \in i\}$ gives
a dense subset.  By relativizing the quantifiers in $A$ to the indices
of non-empty levels, we in effect relativize to a finite subset of the
domain.  We make this more precise:

Suppose $A$ is classically false in some finite structure~\I.  W.l.o.g. we may
assume that the domain of this structure is the naturals $0$, \dots, $n$.  We
extend $\I$ to a $\gdl V$-interpretation $\I^g$ with domain $\N$ as follows:
Since $V$ contains infinitely many values, we can choose $c_1$, $c_2$, $L$ and
$P$ so that $\exists x(x \in i)$ is true for $i = 0$, \dots, $n$ and false
otherwise, and so that $\I^g(\exists x\, P(x)) < 1$.  The number-theoretic
symbols receive their natural interpretation.  The antecedent of $A^g$ clearly
receives the value~$1$, and the consequent receives $\I_g(\exists x\, P(x)) <
1$, so $\I^g \nmodels A^g$.

Now suppose that $\I \nmodels A^g$.  Then $\I(\exists x\, P(x)) < 1$.  In this
case, $\I(x \prec y) = 1$ iff $\I(P(x)) < \I(P(y))$, so $\prec$ defines a
strict order on the domain of $\I$.  It is easily seen that in order for the
value of the antecedent of $A^g$ under \I{} to be greater than that of the
consequent, it must be $= 1$ (the values of all subformulas are either $\le
\I(\exists x\, P(x))$ or $= 1$).  For this to happen, of course, what the
antecedent is intended to express must actually be true in~$\I$, i.e., that $x
\in i$ defines a series of levels and any level $i>0$ is either empty, or for
all~$x$, and~$y$ occuring in some smaller level there is a~$z$ with $x \prec z
\prec y$ and $z \in i$.

To see this, consider the relevant part of the antecedent, $B = \forall
i\bigl[ \forall x, y\forall j\forall k\exists z\,D \lor \forall x\neg(x \in
i)\bigr]$.  If $\I(B) = 1$, then for all $i$, either \( \I(\forall x, y\forall
j\forall k\exists z\,D) = 1 \) or $\I(\forall x\neg(x \in i)) = 1$. In the
first case, we have \( \I(\exists z\,D) = 1 \) for all $x$, $y$, $j$, and $k$.
Now suppose that for all $z$, $\I(D) < 1$, yet $\I(\exists z\, D) = 1$.  Then
for at least some $z$ the value of that formula would have to be $> \I(\exists
z\, P(z))$, which is impossible.  Thus, for every $x$, $y$, $j$, $k$, there is
a $z$ such that $\I(D) = 1$.  But this means that for all $x$, $y$ s.t. $x \in
j$, $y \in k$ with $j, k \le i$ and $x \prec y$ there is a~$z$ with $x \prec z
\prec y$ and $z \in i+1$.

In the second case, where $\I(\forall x\neg(x \in i)) = 1$, we have
that $\I(\neg (x \in i)) = 1$ for all $x$, hence $\I(x \in i) = 0$ and
level~$i$ is empty.

Note that the non empty levels can be distributed over the whole range
of the non-standard model, but since $V$ contains no dense subset, the
total number of non empty levels is finite. Thus, $A$ is false in the
classical interpretation $\I^c$ obtained from $\I$ by restricting $\I$
to the domain $\{ i \suchthat \exists x (x\in i)\}$ and 
$\I^c(Q) = \I(\neg\neg Q)$ for atomic $Q$.
\end{proof}

This shows that no infinite-valued Gödel logic whose set of truth values does
not contain a dense subset, i.e., no countably infinite Gödel logic is
axiomatizable.  We strengthen this result in Section~\ref{sec:prenex} to show
that the prenex fragments are likewise not axiomatizable.


\section{Uncountable Gödel sets}

\subsection{$0$ is contained in the perfect kernel\label{ssec:perfect}}

If $V$ is uncountable, and $0$ is contained in $V^\infty$, then $\gdl V$ is
axiomatizable.  Indeed, Theorem~\ref{thm:gs-id} showed that the sets of
validities of all such~$V$ coincide.  Thus, it is only necessary to establish
completeness of the axioms system $\H$ with respect to $\gdl \bbR$. This
result has been shown by several people over the years.  We give here a
generalization of the proof of Takano \cite{takano}.

\begin{theorem}[Strong completeness of Gödel logic \cite{takano}]
  \label{thm:takano}\hskip0pt plus 3pt minus0pt
  If $\Gamma \entails A$ in $\gdl \R$, then 
  $\Gamma \proves_\H A$.
\end{theorem}

\begin{proof}
  Assume that $\Gamma \nproves A$, we construct an interpretation~$\I$ in
  which $\I(A) = 1$ for all $B \in \Gamma$ and $\I(A) < 1$. Let $y_1$, $y_2$,
  \dots{} be a sequence of free variables which do not occur in $\Gamma \cup
  \Delta$, let $\cT$ be the set of all terms in the language of $\Gamma \cup
  \Delta$ together with the new variables~$y_1$, $y_2$, \dots, and let $\cF =
  \{F_1, F_2, \ldots\}$ be an enumeration of the formulas in this language in
  which $y_i$ does not appear in $F_1$, \dots, $F_i$ and in which each formula
  appears infinitely often.
  
  If $\Delta$ is a set of formulas, we write $\Gamma \Proves \Delta$ if for
  some $A_1$, \dots, $A_n \in \Gamma$, and some $B_1$, \dots, $B_m \in
  \Delta$, $\proves_\H (A_1 \land \ldots \land A_n) \limp (B_1 \lor \ldots
  \lor B_m)$ (and $\nProves$ if this is not the case).  We define a sequence
  of sets of formulas $\Gamma_n$, $\Delta_n$ such that $\Gamma_n \nProves
  \Delta_n$ by induction. First, $\Gamma_0 = \Gamma$ and $\Delta_0 = \{A\}$.
  By the assumption of the theorem, $\Gamma_0 \nProves \Delta_0$.
  
  If $\Gamma_n \Proves \Delta_n \cup \{F_n\}$, then $\Gamma_{n+1} = \Gamma_n
  \cup \{F_n\}$ and $\Delta_{n+1} = \Delta_n$.  In this case, $\Gamma_{n+1}
  \nProves \Delta_{n+1}$, since otherwise we would have $\Gamma_n \Proves
  \Delta_n \cup \{F_n\}$ and $\Gamma_{n} \cup \{F_n\} \Proves \Delta_n$.  But
  then, we'd have 
  that $\Gamma_n \Proves \Delta_n$, which contradicts the induction hypothesis
  (note that $\proves_\H (A \limp B \lor F) \limp ((A \land F \limp B) \limp
  (A \limp B))$).
 
  If $\Gamma_n \nProves \Delta_n \cup \{F_n\}$, then $\Gamma_{n+1} = \Gamma_n$
  and $\Delta_{n+1} = \Delta_n \cup \{F_n, B(y_n)\}$ if $F_n \equiv \forall
  x\, B(x)$, and $\Delta_{n+1} = \Delta_n \cup \{F_n\}$ otherwise.  In the
  latter case, it is obvious that $\Gamma_{n+1} \nProves \Delta_{n+1}$.  In the
  former, observe that by $I10$ and $\axqs$, if $\Gamma_{n} \Proves \Delta_n
  \cup \{\forall x\, B(x), B(y_n)\}$ then also $\Gamma_n \Proves \Delta_n \cup
  \{\forall x\, B(x)\}$ (note that $y_n$ does not occur in $\Gamma_n$ or
  $\Delta_n$).
  
  Let $\Gamma^* = \bigcup_{i=0}^\infty \Gamma_i$ and $\Delta^* =
  \bigcup_{i=0}^\infty \Delta_i$. We have:
  \begin{compactenum}
  \item $\Gamma^* \nProves \Delta^*$, for otherwise there would be a $k$ so
    that $\Gamma_k \Proves \Delta_k$.
  \item $\Gamma \subseteq \Gamma^*$ and $\Delta \subseteq \Delta^*$ (by
    construction).
  \item $\Gamma^* = \cF \setminus \Delta^*$, since each $F_n$ is either in
    $\Gamma_{n+1}$ or $\Delta_{n+1}$, and if for some $n$, $F_n \in \Gamma^*
    \cap \Delta^*$, there would be a~$k$ so that $F_n \in \Gamma_k \cap
    \Delta_k$, which is impossible since $\Gamma_k \nProves \Delta_k$.
  \item If $\Gamma^* \Proves B_1 \lor \ldots \lor B_n$, then~$B_i \in \Gamma^*$
    for some~$i$.  For suppose not, then for $i = 1$, \dots, $n$, $B_i
    \notin \Gamma^*$, and hence, by (3), $B_i \in \Delta^*$.  But then
    $\Gamma^* \Proves \Delta^*$, contradicting~(1).
  \item If $B(t)\in \Gamma^*$ for every~$t \in \cT$, then~$\forall x\, B(x)
    \in \Gamma^*$.  Otherwise, by (3), $\forall x\, B(x) \in \Delta^*$ and so
    there is some $n$ so that $\forall x\, B(x) = F_n$ and $\Delta_{n+1}$
    contains $\forall x\, B(x)$ and $B(y_n)$.  But, again by (3), then $B(y_n)
    \notin \Gamma^*$.
  \item $\Gamma^*$ is closed under provable implication, since if $\Gamma^*
    \Proves A$, then $A \notin \Delta^*$ and so, again by (3), $A \in
    \Gamma^*$.  In particular, if $\proves_\H A$, then $A \in \Gamma^*$.
  \end{compactenum}
  
  Define relations~$\lee$ and~$\equiv$ on~$\cF$ by
  \[
  B\lee C \dequiv B \limp C \in \Gamma^* \quad\text{and}\quad 
  B\equiv C \dequiv B \lee C\land C \lee B.
  \] 
  Then~$\lee$ is reflexive and transitive, since for every~$B$, $\proves_\H B
  \limp B$ and so $B \limp B \in \Gamma^*$, and if $B \limp C \in \Gamma^*$ and
  $C\limp D \in \Gamma^*$ then $B\limp D \in \Gamma^*$, since $B \limp C, C
  \limp D \Proves B \limp D$ (recall (6) above). Hence, $\equiv$ is an
  equivalence relation on~$\cF$. For every~$B$ in~$\cF$ we let~$\eval{B}$ be
  the equivalence class under~$\equiv$ to which~$B$ belongs, and~$\Feq$ the
  set of all equivalence classes. Next we define the relation~$\le$ on~$\Feq$
  by
  \[
  \eval{B} \le \eval{C} \dequiv B \lee C \dequiv B \limp C \in \Gamma^*.
  \]
  Obviously, $\le$ is independent of the choice of representatives $A$, $B$.

\begin{lemma}\label{lm:feq}
  $\lF$ is a countably linearly ordered structure with
  distinct maximal element~$\eval{\top}$ and minimal
  element~$\eval{\bot}$.
\end{lemma}

\begin{proof} 
  Since~$\cF$ is countably infinite, $\Feq$ is countable. For every~$B$
  and~$C$, $\proves_\H (B \limp C) \lor (C \limp B)$ by $\axlin$, and so
  either~$B \limp C\in \Gamma^*$ or~$C \limp B \in \Gamma^*$ (by (4)),
  hence~$\le$ is linear. For every~$B$, $\proves_\H B \limp \top$ and
  $\proves_\H \bot \limp B$, and so~$B \limp \top \in \Gamma^*$ and $\bot
  \limp B \in \Gamma^*$, hence~$\eval{\top}$ and~$\eval{\bot}$ are the maximal
  and minimal elements, respectively. Pick any $A$ in $\Delta^*$. Since~$\top
  \limp \bot \Proves A$, and~$A \notin \Gamma^*$, $\top \limp
  \bot \notin \Gamma^*$, so~$\eval{\top} \neq \eval{\bot}$.
\end{proof}

We abbreviate~$\eval{\top}$ by~$\one$ and~$\eval{\bot}$ by~$\zero$.

\begin{lemma}\label{lm:propertiesfeq}
  The following properties hold in~$\lF$:
  \begin{compactenum}
  \item $\eval{B} = \one \dequiv B\in \Gamma^*$.
  \item $\eval{B\land C} = \min\{\eval{B},\eval{C}\}$.
  \item $\eval{B\lor C} = \max\{\eval{B},\eval{C}\}$.
  \item $\eval{B\limp C} = \one$ if $\eval{B} \le \eval{C}$, $\eval{B \limp
      C}=\eval{C}$ otherwise.
  \item $\eval{\lnot B} = \one$ if $\eval{B} = \zero$; $\eval{\lnot B} =
    \zero$ otherwise.
  \item $\eval{\qe xB(x)} = \sup\{\eval{B(t)} \suchthat t\in\cT\}$.
  \item $\eval{\qa xB(x)} = \inf\{\eval{B(t)} \suchthat t\in\cT\}$.

  \end{compactenum}
\end{lemma}

\begin{proof} 
  (1) If $\eval{B} = \one$, then $\top \impl B \in \Gamma^*$, and hence $B \in
  \Gamma^*$.  And if $B \in \Gamma^*$, then $\top \limp B \in \Gamma^*$ since
  $B \Proves \top \impl B$.  So $\eval{\top} \le \eval{B}$.  It follows that
  $\eval{\top} = \eval{B}$ as also $\eval{B} \le \eval{\top}$.
  
  (2) From~$\Proves B \land C \limp B$, $\Proves B \land C \limp C$ and $D
  \limp B, D \limp C \Proves D \limp B\land C$ for every~$D$, it follows that
  $\eval{B \land C} = \inf \{\eval{B}, \eval{C}\}$, from which (2) follows
  since~$\le$ is linear. (3) is proved analogously.
  
  (4) If $\eval{B} \le \eval{C}$, then $B \limp C \in \Gamma^*$, and since
  $\top \in \Gamma^*$ as well, $\eval{B \limp C} = \one$.  Now suppose that
  $\eval{B} \nleq \eval{C}$. From $B \land (B \limp C) \Proves C$ it follows
  that $\min \{\eval{B}, \eval{B \limp C}\} \le \eval{C}$.  Because $\eval{B}
  \nleq \eval{C}$, $\min \{\eval{B}, \eval{B \limp C}\} \neq \eval{B}$, hence
  $\eval{B \limp C} \le \eval{C}$.  On the other hand, $\proves C \limp (B
  \limp C)$, so $\eval{C} \le \eval{B \limp C}$.
  
  (5) If $\eval{B} = \zero$, $\lnot B = B \impl \bot \in \Gamma^*$, and hence
  $\eval{\lnot B} = \one$ by (1).  Otherwise, $\eval{B} \nleq \eval{\bot}$,
  and so by (4), $\eval{\lnot B} = \eval{B \limp \bot} = \zero$.
  
  (6) Since $\proves_\H B(t)\limp \qe x\, B(x)$, $\eval{B(t)} \le \eval{\qe
    x\, B(x)}$ for every~$t\in\cT$. On the other hand, for every~$D$ without
  $x$ free,
  \begin{align*}
    &&&\eval{B(t)} \le \eval{D} && \qquad\mbox{for every $t\in\cT$}\\
    &\dequiv&& B(t)\limp D\in \Gamma^*&&\qquad\mbox{for every $t\in\cT$}\\
    &\Rightarrow&& \qa x(B(x)\limp D)\in\Gamma^* &
              & \qquad\text{by property (5) of $\Gamma^*$}\\
    &\Rightarrow&& \qe x\,B(x) \limp D \in\Gamma^* &&\qquad
    \mbox{since $\qa x(B(x)\limp D) \Proves \qe x\,B(x)\limp D$}\\
    &\dequiv&& \eval{\qe x\, B(x)}\le \eval{D}.
  \end{align*}
  (7) is proved analogously.
\end{proof}

$\lF$ is countable, let $\zero = a_0, \one = a_1, a_2, \ldots$ be an
enumeration. Define $h(\zero) = 0$, $h(\one) = 1$, and define~$h(a_n)$
inductively for $n>1$: Let
$ a_n^- = \max\{a_i\suchthat i<n \text{ and } a_i < a_n\}$ and
$ a_n^+ = \min\{a_i\suchthat i<n \text{ and } a_i > a_n\}$, and define 
$h(a_n) = (h(a_n^-) + h(a_n^+))/2$ (thus, $a_2^-=\zero$ and
$a_2^+=\one$ as $\zero=a_0<a_2<a_1=\one$, hence $h(a_2)=\frac12$).
Then $h\colon \lF \to \bbQ \cap [0,1]$ is a strictly monotone map
which preserves infs and sups.  By
Lemma~\ref{lemma:countable-into-cantor} there exists a
\g-embedding~$h'$ from~$\bbQ \cap [0, 1]$ into~$\lR$ which is also
strictly monotone and preserves infs and sups. Put~$\I(B) =
h'(h(\eval{B}))$ for every atomic~$B\in\cF$ and we obtain a
$V_\R$-interpretation.

Note that for every~$B$, $\I(B) = 1$ iff $\eval{B} = \one$ iff $B \in
\Gamma^*$. Hence, we have $\I(B) = 1$ for all $B \in \Gamma$ while if $A
\notin \Gamma^*$, then~$\I(A) < 1$, so~$\Gamma \nentails A$.  Thus we have
proven that on the assumption that if $\Gamma \nproves A$, then $\Gamma
\nentails A$
\end{proof}

As already mentioned we obtain from this completeness proof together with the
soundness theorem (Theorem~\ref{thm:soundness}) and Theorem~\ref{thm:gs-id}
the characterization of recursive axiomatizability:

\begin{theorem}\label{thm:fo:zeroinP}
  Let $V$ be a Gödel set with $0$ contained in the perfect kernel of $V$.
  Suppose that $\Gamma$ is a set of closed formulas. Then $\Gamma \entails_V
  A$ iff $\Gamma \proves_\H A$.
\end{theorem}

\begin{corollary}[Deduction theorem for Gödel logics]\label{thm:deduction}
  Suppose that $\Gamma$ is a set of formulas, and $A$ is a closed formula.
  Then
  \[
  \Gamma, A \proves_\H B \quad \text{iff} \quad \Gamma \proves_\H  A \limp B.
  \]
\end{corollary}

\begin{proof} 
  Use the soundness theorem (Theorem~\ref{thm:soundness}), completeness
  theorem (Theorem~\ref{thm:fo:zeroinP}) and the semantic deduction
  theorem~\ref{thm:semdeduc}.  Another proof would be by induction on the
  length of the proof. See \cite{hajek}, Theorem~2.2.18.
\end{proof}

\subsection{$0$ is isolated\label{ssec:isolated}}

In the case where $0$ is isolated, and thus also not contained in the
perfect kernel, we will transform a counter example in $\GR$ for
$\Gamma,\Pi\entails A$, where $\Pi$ is a set of sentences stating that 
every infimum is a minimum, into a counter example in $\gdl V$ 
for $\Gamma\entails A$.

\begin{lemma}\label{lm:equivaxioms}
  Let $x,\bar y$ be the free variables in $A$.
  \[ 
     \proves_{\Ho} \forall\bar y(\lnot \forall x\, A(x,\bar y) \limp
     \exists x\, \lnot A(x,\bar y))
   \]
\end{lemma}

\begin{proof}
  It is easy to see that in all Gödel logics the following weak form of the
  law of excluded middle is valid: $\lnot\lnot A(a) \lor \lnot A(a)$.  By
  quantification we obtain $\qa x \lnot\lnot A(x) \lor \qe x \lnot A(x)$ and
  by valid quantifier shifting rules $\lnot\lnot \qa x\, A(x) \lor \qe \lnot
  A(x)$.  From the intuitionistically valid $\lnot A\lor B \limp (A \limp B)$
  we can prove $\lnot \qa x\, A(x) \limp \qe x \lnot A(x)$. A final
  quantification of the free variables concludes the proof.
\end{proof}

\begin{theorem}\label{thm:fo:iso}
  Let $V$ be an uncountable Gödel set where $0$ is isolated.  Suppose $\Gamma$
  is a set of closed formulas.  Then $\Gamma \entails_V A$ iff $\Gamma
  \proves_{\Ho} A$.
\end{theorem}

\begin{proof}
  If: Follows from soundness (Theorem~\ref{thm:soundness}) and the observation
  that $\axiso$ is valid for any $V$ where $0$ is isolated.
  
  Only if: We already know from Theorem~\ref{thm:gs-id} that the entailment
  relation of $V$ and $V \cup [\inf P,1]$ coincide, where $P$ is the perfect
  kernel of $V$. So we may assume wthout loss of generality that $V$ already
  is of this form, i.e.\ that $\lambda = \inf P$ and $V \cap [\lambda, 1] =
  [\lambda, 1]$. Let $V' = [0,1]$.  Define
  \[
  \Pi = \{\qa{\bar y}(\lnot\qa x\, A(x,\bar y) \limp \qe x\lnot A(x,\bar y))
  \suchthat A(x,\bar y) \mbox{ formula}\}
  \] 
  where $A(x,\bar y)$ ranges over \emph{all} formulas with free variables $x$
  and $\bar y$. We consider the entailment relation in $V'$. Either $\Pi,
  \Gamma \entails_{V'} A$ or $\Pi, \Gamma \nentails_{V'} A$.  In the former
  case we know from the strong completeness of $\H$ for $\GR$ that there are
  finite subsets $\Pi'$ and $\Gamma'$ of $\Pi$ and $\Gamma$, respectively,
  such that $\Pi', \Gamma' \proves_{\H} A$.  Since all the sentences in $\Pi$
  are provable in $\Ho$ (see Lemma~\ref{lm:equivaxioms}) we obtain that
  $\Gamma' \proves_{\Ho} A$.  In the latter case there is an interpretation
  $\I'$ such that
  \[
   \inf \{ \I'(G) \suchthat G\in \Pi \cup \Gamma \} > \I'(A).
  \]
  It is obvious from the structure of the formulas in $\Pi$ that their truth
  value will always be either $0$ or $1$. Combined with the above we know that
  for all $G \in \Pi$, $\I'(G) = 1$.  Next we define a function $f(x)$ which
  maps values from $\subval(\I', \Gamma \cup \Pi \cup\{A\})$ into $V$:
  \[
     f(x) = \begin{cases} 0 & x=0\\
                          \lambda + x/(1-\lambda) & x>0 
             \end{cases}
  \]
  We see that $f$ satisfies conditions~(1) and~(2) of
  Lemma~\ref{lem:induced-interpretation}, but we cannot use
  Lemma~\ref{lem:induced-interpretation} directly, as not all existing infima
  and suprema are necessarily preserved.
  
  Consider as in Lemma~\ref{lem:induced-interpretation} the interpretation
  $\I_f(B) = f(\I'(B))$ for atomic subformulas of $\Gamma \cup \Pi \cup\{A\}$.
  We want to show that the identity $\I_f(B) = f(\I'(B))$ extends to all
  subformulas of $\Gamma\cup \Pi \cup \{A\}$. For propositional connectives
  and the existentially quantified formulas this is obvious.  The important
  case is $\qa x\, A(x)$. First assume that $\I'(\qa x\, A(x)) > 0$. Then it
  is obvious that $\I_f(\qa x\, A(x)) = f(\I'(\qa x\, A(x)))$. In the case
  where $\I'(\qa x\, A(x)) = 0$ we observe that $A(x)$ contains a free
  variable and therefore $\lnot \qa x\, A(x) \limp \qe x \lnot A(x) \in \Pi$,
  thus $\I'(\lnot \qa x\, A(x) \limp \qe x\lnot A(x)) = 1$. This implies that
  there is a witness $c$ such that $\I'(A(c)) = 0$. Using the induction
  hypothesis we know that $\I_f(A(c)) = 0$, too. We obtain that $\I_f(\qa x\,
  A(x)) = 0$, concluding the proof.
  
  Thus we have shown that $\I_f$ is a counterexample to $\Gamma \entails_V A$
  which completes the proof of the theorem.
\end{proof}

\subsection{$0$ not isolated but not in the perfect kernel}\label{sec:uncount-nonax}

In the preceding sections, we gave axiomatizations for the logics based on
those uncountably infinite Gödel sets~$V$ where $0$ is either isolated or in
the perfect kernel of~$V$.  It remains to determine whether logics based on
uncountable Gödel sets where $0$ is neither isolated nor in the perfect kernel
are axiomatizable.  The answer in this case is negative.  If $0$ is not
isolated in $V$, $0$ has a countably infinite neighborhood.  Furthermore, any
sequence $(a_n)_{n \in \N} \to 0$ is so that, for sufficiently large $n$, $V
\cap [0, a_n]$ is countable and hence, by (the proof of)
Theorem~\ref{bendixon}, contains no densely ordered subset.  This fact is the
basis for the following non-axiomatizability proof, which is a variation on
the proof of Theorem~\ref{thm:count-nonax}.

\begin{theorem}\label{thm:uncount-noniso-nonax}
  If $V$ is uncountable, $0$ is not isolated in $V$, but not in the perfect
  kernel of~$V$, then $\gdl{V}$ is not axiomatizable.
\end{theorem}

\begin{proof}
  We show that for every sentence $A$ there is a sentence $A^h$ s.t. $A^h$ is
  valid in $\gdl V$ iff $A$ is true in every finite (classical) first-order
  structure.
  
  The definition of $A^h$ mirrors the definition of $A^g$ in the proof of
  Theorem~\ref{thm:count-nonax}, except that the construction there is carried
  out infinitely many times for $V \cap [0, a_n]$, where $(a_n)_{n \in \N}$ is
  a strictly descending sequence, $a_n > 0$ for all $n$, which converges
  to~$0$.  Let $P$ be a binary and $L$ be a ternary predicate symbol not
  occurring in~$A$ and let $R_1$, \dots, $R_n$ be all the predicate symbols
  in~$A$.  We use the abbreviations $x \in_\ell y \equiv \neg\neg L(x, y,
  \ell)$ and $x \prec_\ell y \equiv (P(y, \ell) \impl P(x, \ell)) \impl P(y,
  \ell)$.  As before, for a fixed $\ell$, provided $\I(\exists x\, P(x, \ell))
  < 1$, $\I(x \prec_\ell y) = 1$ iff $\I(P(x, \ell)) < \I(P(y, \ell))$, and
  $\I(x \in_\ell y)$ is always either $0$ or $1$.  We also need a binary
  predicate symbol $Q(\ell)$ to give us the descending
  sequence~$(a_n)_{n\in\N}$: Note that $\I(\neg \forall \ell\, Q(\ell)) = 1$
  iff $\inf \{\I(Q(d)) : d \in \card{\I}\} = 0$ and $\I(\exists \ell\, \neg
  Q(\ell)) = 1$ iff $0 \notin \{\I(Q(d)) : d \in \card{\I}\}$.
  
  Let $A^h \equiv$
  \begin{equation}\label{fm:0notiso}
    \left\{
      \begin{array}{l}
        S \land  
        \forall \ell((Q(s(\ell)) \impl Q(\ell)) \impl Q(s(\ell)) \land {}\\
        \quad \neg \forall \ell\, Q(\ell) \land \exists \ell\, \neg Q(\ell) \land {}\\
        \quad \forall \ell\forall x((Q(\ell) \impl P(x, \ell)) \impl Q(\ell)) \land {}\\ 
        \quad \forall \ell \exists x\exists y(x \in_\ell 0 \land {} y \in_\ell 0 \land 
        x \prec_\ell y) \land {}\\
        \quad 
        \quad \forall \ell\forall i\bigl[\forall x, y \forall j \forall k \exists z\, E
        \lor \forall x\neg(x \in_\ell s(i))\bigr]
      \end{array}
    \right\}
    \impl (A' \lor \exists \ell \exists u\, P(u, \ell) \lor \exists \ell\, Q(\ell))
  \end{equation}
  where $S$ is the conjunction of the standard axioms for $0$, successor and
  $\le$, with double negations in front of atomic formulas,
  \[ 
  E\equiv \begin{array}{l}
    (j \le i \land x \in_\ell j \land k\le i \land y \in_\ell k \land x \prec_\ell y) \impl {}\\
    \qquad \impl (z \in_\ell s(i) \land x \prec_\ell z \land z \prec_\ell y)
  \end{array}
  \]
  and $A'$ is $A$ where every atomic formula is replaced by its double
  negation, and all quantifiers are relativized to the predicate
  $R(\ell)\equiv \forall i\exists x(x \in_\ell i)$.
  
  The idea here is that an interpretation $\I$ will define a sequence
  $(a_n)_{n\in\N} \to 0$ by $a_n = \I(Q(\bar{n}))$ where $a_n > a_{n+1}$, and
  $0 < a_n < 1$ for all~$n$.  Let $L_\ell^i = \{x : \I(x \in_\ell i)\}$ be the
  $i$-th $\ell$-level.  $P(x, \ell)$ orders the set $\bigcup_i L_\ell^i = \{x
  : \I(\exists i\, x \in_\ell i) = 1\}$ in a subordering of $V \cap [0, a_n]$:
  $x \prec_\ell y$ iff $\I(x \prec_\ell y) = 1$.  Again we force that whenever
  $x, y \in L_\ell^i$ with $x \prec_\ell y$, there is a $z \in L_\ell^{i+1}$
  with $x \prec_\ell z \prec_\ell y$, or, if no possible such $z$ exists,
  $L_\ell^{i+1} = \emptyset$.  Let $r(\ell)$ be the least $i$ so that
  $L_\ell^i$ is empty, or $\infty$ otherwise.  If $r(\ell) = \infty$ then
  there is a densely ordered subset of $V \cap [0, a_\ell]$.  So if $0$ is not
  in the perfect kernel, for some sufficiently large $L$, $r(\ell) < \infty$
  for all $\ell > L$.  $\I(R(\ell)) = 1$ iff $r(\ell) = \infty$ hence $\{\ell
  : \I(R(\ell)) = 1\}$ is finite whenever the interpretations of $P$, $L$, and
  $Q$ are as intended.
  
  Now if $A$ is classically false in some finite structure~\I, we can again
  choose a $\gdl{V}$-interpretation $\I^h$ in which the interpretations of
  $P$, $Q$, $L$ are as intended, the number theoretic predicates and functions
  receive their standard interpretation, there are as many $\ell$ with
  $\I^h(R(\ell)) = 1$ as there are elements in the domain of $\I$, and the
  predicates of $A$ behave on $\{\ell : \I(R(\ell)) = 1\}$ just as they do
  on~$\I$.  $\I^h \nmodels A^h$.
  
  On the other hand, if $\I \nmodels A^h$, then the value of the consequent is
  $< 1$.  Then as required, for all $x$, $\ell$, $\I(P(x, \ell)) < 1$ and
  $\I(Q(\ell)) < 1$.  Since the antecedent, as before, must be $=1$, this
  means that $x \prec_\ell y$ expresses a strict ordering of the elements of
  $L_\ell^i$ and $\I(((Q(s(\ell) \impl Q(\ell)) \impl Q(s(\ell))) = 1$ for all
  $\ell$ guarantees that $\I(Q(s(\ell))) = a_{n+1} < a_n = \I(Q(\ell))$.  The
  other conditions are likewise seen to hold as intended, so that we can
  extract a finite countermodel for $A$ based on the interpretation of the
  predicate symbols of $A$ on $\{\ell : \I(R(\ell)) = 1\}$, which must be
  finite.
\end{proof}

\section{Fragments}\label{sec:fragments}

\fxnote{Für welche Fragmente kann/soll man starke vollständigkeit (im
  axiomatisierbaren Fall) zeigen?}

\subsection{Prenex fragments}\label{sec:prenex}

One interesting restriction of the axiomatizability problem is the question
whether the prenex fragment of $\gdl{V}$, i.e., the set of prenex formulas
valid in $\gdl{V}$, is axiomatizable.  This is non-trivial, since in general
in Gödel logics, arbitrary formulas are not equivalent to prenex formulas.
Thus, so far the proofs of non-axiomatizability of the logics treated in
Sections~\ref{sec:countable} and~\ref{sec:uncount-nonax} do not establish the
non-axiomatizability of their prenex fragments, nor do they exclude the
possibility that the corresponding prenex fragments are r.e.  We investigate
this question in this section, and show that the prenex fragments of all
finite and uncountable Gödel logics are r.e., and that the prenex fragments of
all countably infinite Gödel logics are not r.e.  The axiomatizability result
is obtained from a version of Herbrand's Theorem for finite and
uncountably-valued Gödel logics, which is of independent interest.  The
non-axiomatizability of countably infinite Gödel logics is obtained as a
corrolary of Theorem~\ref{thm:count-nonax}.

Let $V$ be a Gödel set which is either finite or uncountable. Let $\Gd$ be a
Gödel logic with such a truth value set. We show how to effectively associate
with each prenex formula $A$ a quantifier-free formula $A^\ast$ which is valid
in $\Gd$ if and only if $A$ is a tautology.  The axiomatizability of the
prenex fragment of $\Gd$ then follows from the axiomatizability of \LC{} (in
the infinite-valued case) and propositional $\gdl m$ (in the finite-valued
case).

\begin{definition}[Herbrand form]
Given a prenex formula $A \equiv \Q_1
x_1\ldots\Q_n x_n\,B(\bar x)$ ($B$ quantifier free), the
\emph{Herbrand form}~$A^H$ of $A$ is $\exists x_{i_1}\ldots
\exists x_{i_m}\, B(t_1, \ldots, t_n)$, where $\{x_{i_j} : 1 \le j \le
m\}$ is the set of existentially quantified variables in $A$, and
$t_i$ is $x_{i_j}$ if $i = i_j$, or is $f_i(x_{i_1}, \ldots, x_{i_k})$
if $x_i$ is universally quantified and $k = \max \{j : i_j < i\}$. We
will write $B(t_1,\ldots,t_n)$ as $B^F(x_{i_1}, \ldots, x_{i_m})$ if
we want to emphasize the free variables.
\end{definition}

\begin{lemma}
\label{skolemization}
If $A$ is prenex and $\modelsg A$, then $\modelsg A^H.$
\end{lemma}
 
\begin{proof}
Follows from the usual laws of quantification, which are valid in all Gödel
logics.
\end{proof}

Our next main result will be Herbrand's theorem for $\gdl V$ for $V$
uncountable or finite.  The {\em Herbrand universe} $\HU(B^F)$ of
$B^F$ is the set of all variable-free terms which can be constructed
from the set of function symbols occurring in $B^F$. To prevent
$\HU(B^F)$ from being finite or empty we add a constant and a function
symbol of positive arity if no such symbols appear in~$B^F$.  The {\em
  Herbrand base} $\HB(B^F)$ is the set of atoms constructed from the
predicate symbols in $B^F$ and the terms of the Herbrand universe. In
the next theorem we will consider the Herbrand universe of a formula
$\exists \tup{x}\, B^F(\tup{x})$.  We fix a non-repetitive enumeration
$C_1$, $C_2$, \dots of $\HB(B^F)$, and let 
$X_\ell = \{\bot, C_1, \ldots, C_\ell, \top\}$ (we may take $\top$ to
be a formula which is always $=1$). $B^F(\tup{t})$ is an
\emph{$\ell$-instance} of $B^F(\tup{x})$ if the atomic subformulas of
$B^F(\tup{t})$ are in $X_\ell$. 

\begin{definition} 
  An $\ell$-\emph{constraint} is a non-strict linear ordering
  $\preceq$ of $X_\ell$ s.t. $\bot$ is minimal and $\top$ is maximal.
  An interpretation~$\I$ {\em fulfils} the constraint $\preceq$
  provided for all $C, C' \in X_\ell$, $C \preceq C'$ iff 
  $\I(C) \le \I(C')$.  We say that the constraint $\preceq'$ on
  $X_{\ell+1}$ \emph{extends} $\preceq$ if for all $C, C' \in X_\ell$,
  $C \preceq C'$ iff $C \preceq' C'$.
\end{definition}

Lemma~\ref{lem:g-embed} showed that if $h: V \to W$ is a \g-embedding and $\I$
is a $V$-interpretation, then $h(\I(A)) = \I_h(A)$ for any formula~$A$.  If no
quantifiers are involved in $A$, this also holds without the requirement of
continuity.  For the following proof we need a similar notion.  Let $V$ be a
Gödel sets, $X$ a set of atomic formulas, and suppose there is an
order-preserving, strictly monotone $h\colon \{\I(C)\colon C \in X\} \to V$
which is so that $h(1) = 1$ and $h(0) = 0$.  Call any such $h$ a \emph{truth
value injection on $X$}.  Now suppose $B$ is a quantifier-free formula, and
$X$ its set of atomic subformulas.  Two interpretations $\I$, $\J$ are
\emph{compatible on $X$} if $\I(C) \le \I(C')$ iff $\J(C) \le \J(C')$ for all
$C \in X$.

\begin{proposition}\label{tviso}
 Let $B^F$ be a quantifier free formula, and $X$ its set of atomic subformulas
 together with $\top$, $\bot$.  If $\I$, $\J$ are compatible on $X$, then
 there is a truth value injection~$h$ on $X$ with $h(\I(C^F)) = \J(C^F)$.
\end{proposition}

\begin{proof}
Let $h(\I(C)) = \J(C)$ for $B \in X$.  Since $\I$, $\J$ are compatible on $X$,
$\I(C) \le \I(C')$ iff $\J(C) \le \J(C')$, and hence $\I(C) \le \I(C')$ iff
$h(\I(C)) \le h(\J(C'))$ and $h$ is strictly monotonic.  The conditions $h(0)
= 0$ and $h(1) = 1$ are satisfied by definition, since $\top$, $\bot \in X$.
We get $h(\I(B^F)) = \J(B^F)$ by induction on the complexity of~$A$.
\end{proof}

\begin{proposition}\label{ellI}
(a) If $\preceq'$ extends $\preceq$, then every $\I$ which fulfills $\preceq'$
also fulfills $\preceq$. (b) If $\I$, $\J$ fulfill the $\ell$-constraint
$\preceq$, then there is a truth value injection~$h$ on $X_\ell$ with
$h(\I(B^F(\tup{t}))) = \J(B^F(\tup{t}))$ for all $\ell$-instances
$B^F(\tup{t})$ of $B^F(\tup{x})$; in particular, $\I(B^F(\tup{t})) = 1$ iff
$\J(B^F(\tup{t})) = 1$.
\end{proposition}

\begin{proof}
(a) Obvious. (b) Follows from Proposition~\ref{tviso} together with the
  observation that $\I$ and $\J$ both fulfill $\preceq$ iff they are
  compatible on $X_\ell$.
\end{proof} 

\begin{lemma} \label{univH}
Let $B^F$ be a quantifier-free formula, and let $V$ be a finite or uncountably
infinite Gödel set. If $\modelsg \exists \tup{x}\, B^F(\tup{x})$ then there are
tuples $\tup{t}_1, \dots \overline{t}_n$ of terms in $U(B^F)$, such that
$\modelsg \bigvee_{i=1}^{n} B^F(\tup{t}_i)$.
\end{lemma}

\begin{proof}
Suppose first that $V$ is uncountable.  By Theorem~\ref{bendixon}, $V$
contains a dense linear subordering. We construct a ``semantic tree'' $\T$;
i.e., a systematic representation of all possible order types of
interpretations of the atoms $C_i$ in the Herbrand base.  $\T$ is a rooted
tree whose nodes appear at levels.  Each node at level~$\ell$ is labelled with
an $\ell$-constraint.

$\T$ is constructed in levels as follows: At level~$0$, the root of~$\T$ is
labelled with the constraint $\bot < \top$.  Let $\nu$ be a node added at
level $\ell$ with label $\preceq$, and let $T_\ell$ be the set of terms
occurring in $X_\ell$.  Let (*) be: For every interpretation~$\I$ which
fulfils~$\preceq$, there is some $\ell$-instance $B^F(\tup{t})$ so that
$\I(B^F(\tup{t}))=1$. If (*) obtains, $\nu$ is a leaf node of~$\T$, and no
successor nodes are added at level $\ell +1$.
  
Note that by Proposition~\ref{ellI}(b), any two
interpretations which fulfill $\preceq$ make the same $\ell$-instances
of $B^F(\tup{t})$ true; hence $\nu$ is a leaf node if and only if there
is an $\ell$-instance $A(\tup{t})$ s.t. $\I(A(\tup{t})) = 1$ for
all interpretations $\I$ that fulfil~$\preceq$.

If (*) does not obtain, for each $(\ell + 1)$-constraint $\preceq'$
extending $\preceq$ we add a successor node $\nu'$ labelled with
$\preceq'$ to $\nu$ at level $\ell+1$.

We now have two cases:

(1) $\T$ is finite. Let $\nu_1, \ldots, \nu_m$ be the leaf nodes of~$\T$ of
    levels $\ell_1$, \ldots, $\ell_m$, each labelled with a constraint
    $\preceq_1$, \dots, $\preceq_m$.  By (*), for each $j$ there is an
    $\ell_j$-instance $B^F(\tup{t}_j)$ with $\I(B^F(\tup{t})) = 1$ for all
    $\I$ which fulfill~$\preceq_j$.  It is easy to see that every
    interpretation fulfills at least one of the $\preceq_j$.  Hence, for all
    $\I$, $\I(B^F(\tup{t}_1) \lor \ldots \lor B^F(\tup{t}_m)) = 1$, and so
    $\modelsg \bigvee_{i=1}^m B^F(\tup{t}_i)$.

(2) $\T$ is infinite.  By König's lemma, $\T$ has an infinite branch
    with nodes $\nu_0$, $\nu_1$, $\nu_2$, \dots where $\nu_\ell$ is
    labelled by $\preceq_\ell$ and is of level $\ell$.  Each
    $\preceq_{\ell+1}$ extends $\preceq_\ell$, hence we can form
    $\mathnormal{\preceq} = \bigcup_\ell \mathnormal{\preceq_\ell}$.  
    Let $V' \subseteq V$ be a
    non-trivial densely ordered subset of $V$, let $V' \ni c < 1$, and
    let $V'' = V' \cap [0, c)$.  $V''$ is clearly also densely
    ordered. Now let $V_c$ be $V'' \cup \{0, 1\}$, and let $h: B(A(x))
    \cup \{\bot, \top\} \to V_c$ be an injection which is so that, for
    all $A_i, A_j \in B(A(x))$, $h(A_i) \le h(A_j)$ iff $A_i \preceq
    A_j$, $h(\bot) = 0$ and $h(\top) = 1$. We define an interpretation
    $\I$ by: $f^{\I}(t_1, \ldots, t_n) =
    f(t_1, \ldots, t_n)$ for all $n$-ary function symbols~$f$ and
    $P^{\I}(t_1, \ldots, t_n) = h(P(t_1, \ldots, t_n))$ for all $n$-ary
    predicate symbols~$P$ (clearly then, $\I(A_i) = h(A_i)$). By
    definition, $\I$ $\ell$-fulfills $\preceq_\ell$ for all $\ell$. By
    (*), $\I(A(\tup{t})) < 1$ for all $\ell$-instances $A(\tup{t})$ of
    $A(x)$, and by the definition of $V_c$, $\I(A(\tup{t})) < c$.
    Since every $A(\tup{t})$ with $\tup{t} \in U(A(x))$ is an
    $\ell$-instance of $A(x)$ for some $\ell$, we have $\I(\exists x\,
    A(\tup{x})) \le c < 1$.This contradicts the assumption that
    $\modelsg \exists \tup{x}\, A(\tup{x})$.

If $V$ is finite, the proof is the similar, except simpler.  Suppose
$\card{V} = n$.  Call a constraint $\preceq$
\emph{$n$-admissible} if there is some $V$-interpretation $\I$ which
fulfills it.  Such $\preceq$ have no more than $n$ equivalence classes under
the equivalence relation $C \sim C'$ iff $C \preceq C'$ and $C' \preceq C$.
In the construction of the semantic tree above, replace each mention of
$\ell$-constraints by $n$-admissible $\ell$-constraints.  The argument in the
case where the resulting tree is finite is the same.  If $\T$ is infinite,
then the resulting order ${\preceq} = \bigcup_\ell {\preceq_\ell}$ is
$n$-admissible, since all $\preceq_\ell$ are. Let $c = \max \{ b: b \in V, b <
1\}$ and $V_c = V$.  The rest of the argument goes through without change.
\end{proof}

\begin{lemma}\label{lem:reskolem}
Let $\exists \tup{x}\,
B^F(\tup{x})$ be the Herbrand form of the prenex formula $A \equiv \tup{\Q}_i
B(\tup{y}_i)$, and let $\tup{t}_1, \ldots, \tup{t}_m$ be tuples of terms in
$\HU(B^F)$. If $\gdl{V} \models \bigvee_{i=1}^{m} B^F(\tup{t_i})$, then
$\gdl{V} \models A.$
\end{lemma}

\begin{proof} For any Gödel set $V$, the following rules are valid
  in~$\gdl{V}$:

(1) $A \lor B  \vdash B \lor A$.

(2) $(A \lor B) \lor C \vdash A \lor (B \lor C)$ .

(3) $A \lor (B \lor B) \vdash A \lor B$ .

(4) $A(y) \vdash \forall x\, A(x)$.

(5) $A(t) \vdash \exists x\,  A(x)$.

(6) $\forall x (A(x) \lor B) \vdash \forall x\, A(x) \lor B$.

(7) $\exists x (A(x) \lor B) \vdash \exists x\, A(x) \lor B$.

\noindent ($x$ is not free in $B$.)  The result follows from
\cite{BCF01LPAR}, Lemma~6, and are also easily verified directly.
\end{proof}

\begin{theorem} \label{HT}
Let $A$ be prenex, $\exists \tup{x}\, B^F(\tup{x})$ its Herbrand form, and let
$V$ be a finite or uncountably infinite Gödel set. Then $\modelsg A$ iff there
are tuples $\tup{t}_1, \ldots \tup{t}_m$ of terms in $\HU(B^F)$, such that
$\modelsg \bigvee_{i=1}^{m} B^F(\tup{t}_i)$.
\end{theorem}

\begin{proof}
If: This is Lemma~\ref{lem:reskolem}.  Only if: By Lemma~\ref{skolemization}
and Lemma~\ref{univH}.
\end{proof}

\begin{remark} 
An alternative proof of Herbrand's theorem can be obtained using the
analytic calculus {\it HIF} (``Hypersequent calculus for
Intuitionistic Fuzzy logic'') \cite{BaazZach00CSL}.
\end{remark}

\begin{theorem}
The prenex fragment of a Gödel logic based on a truth value set~$V$ which is
either finite or uncountable infinite is axiomatizable. An axiomatization is
given by the standard axioms and rules for~$\LC$ extended by the rules
(4)--(7) of the proof of Lemma~\ref{lem:reskolem}. For the
$m$-valued case add the characteristic axiom for $\gdl m$, \( G_m \equiv
\bigvee_{i=1}^m\bigvee_{j=i+1}^{m+1} ((A_i \impl A_j) \land (A_j \impl A_i)).
\)
\end{theorem}

\begin{proof}  
Completeness: Let $\tup{\Q} \tup{y}_i B(\tup{y})$ be a prenex formula valid in
$\gdl V$. By Theorem~\ref{HT}, a Herbrand disjunction $\bigvee_{i=1}^{n}
B^F(\tup{t}_i)$ is a tautology in $\gdl{V}$.  Hence, it is provable in \LC{}
or $\LC + G_m$ \cite[Chapter 10.1]{gottwald}.  $\tup{\Q} \tup{y} B(\tup{y})$
is provable by Lemma~\ref{lem:reskolem}.

Soundness: The rules in the proof of Lemma~\ref{lem:reskolem} are valid in
$\gdl{V}$. In particular, note that $\forall x(A(x) \lor B) \to (\forall x\,
A(x) \lor B)$ with $x$ not free in $B$ is valid in all Gödel logics, and
$\exists x(A(x) \lor B) \to \exists x\, A(x) \lor B$ is already
intuitionistically valid.
\end{proof}

In Theorem~\ref{thm:count-nonax}, we showed that for every first-order formula
$A$, there is a formula $A^g$ which is valid in $\gdl{V}$ for $V$ countably
infinite iff $A$ is valid in every finite classical interpretation.  We now
strengthen this result to show that the prenex fragment of $\gdl{V}$ (for $V$
countably infinite) is likewise not axiomatizable.  This is done by showing
that if $A$ is prenex, then there is a formula $A^G$ which is also prenex and
which is valid in $\gdl{V}$ iff $A^g$ is.  Note that not all quantifier
shifting rules are generally valid in Gödel logics, so we have to show that
for the particular case of formulas of the form of $A^g$, there is a prenex
formula which is valid in $\gdl{V}$ iff $A^g$ is.

\begin{theorem}
If $V$ is countably infinite, the prenex fragment of $\gdl V$ is not
r.e.
\end{theorem}

\begin{proof}
By the proof of Theorem~\ref{thm:count-nonax}, a formula $A$ is true in all
finite models iff $\gdl V \models A^g$.  $A^g$ is of the form $B \impl (A' \lor
\exists u\, P(u))$.  We show that $A^g$ is validity-equivalent in $\gdl V$ to a
prenex formula.


From Lemma~\ref{lm:crisp} we see that each crisp formula is equivalent to a
prenex formula; let $A_0$ be a prenex form of $A'$.  Since all quantifier
shifts for conjunctions are valid, the antecedent~$B$ of $A^g$ is equivalent
to a prenex formula $\Q_1x_1 \ldots \Q_n x_n B_0(x_1, \ldots, x_n)$.  Hence,
$A^g$ is equivalent to $\tup{\Q}\tup{x} B_0(\tup{x}) \to (A_0 \lor \exists u\,
P(u))$.

Let $\Q_i'$ be $\exists$ if $\Q_i$ is $\forall$, and $\forall$ if
$\Q_i$ is $\exists$, let $C \equiv A_0 \lor \exists u\, P(u)$, and $v
= \I(\exists u\, P(u)))$.  We show that $\tup{\Q}\tup{x}\,
B_0(\tup{x}) \to C$ is equivalent to $\tup{\Q}'\tup{x}(B_0(\tup{x})
\to C)$ by induction on~$n$.  Let $\tup{\Q}\tup{x}B_0 \equiv \Q_1
x_1\ldots \Q_{i}x_{i} B_1(d_1, \ldots, d_{i-1}, x_i)$.  Since
quantifier shifts for $\exists$ in the antecent of a conditional are
valid, we only have to consider the case $\Q_i = \forall$. Suppose
$\I(\forall x_i\, B_1(\tup{d}, x_i) \impl C) \neq \I(\exists
x_i(B_1(\tup{d}, x_i) \impl C)$. This can only happen if $\I(\forall
x_i\, B_1(\tup{d}, x_i)) = \I(C) < 1$ but $\I(B_1(\tup{d}, c)) > \I(C)
\ge v$ for all $c$.  However, it is easy to see by inspecting $B$ that
$\I(B_1(\tup{d}, c))$ is either $=1$ or $\le v$.

Now we show that $\I(B_0(\tup{d}) \impl (A_0 \lor \exists u\, P(u))) =
\I(\exists u(B_0(\tup{d}) \impl (A_0 \lor P(u))))$.  If $\I(A_0) = 1$,
then both sides equal $= 1$.  If $\I(A_0) = 0$, then $\I(A_0 \lor
\exists u\, P(u)) = v$.  The only case where the two sides might
differ is if $\I(B_0(\tup{d})) = v$ but $\I(A_0 \lor P(c)) = \I(P(c))
< v$ for all $c$.  But inspection of $B_0$ shows that $\I(B_0(\tup{t}))
= 1$ or $= \I(P(e))$ for some $e \in \tup{d}$ (the only subformulas of
$B_0(\tup{d})$ which do not appear negated are of the form $e' \prec
e$).  Hence, if $\I(B_0(\tup{d})) = v$, then for some~$e$, $\I(P(e)) =
v$.

Last we consider the quantifiers in $A_0 \equiv \tup{\Q} \tup{y}\, A_1$.
Since $A_0$ is crisp, $\I(B_0(\tup{d}) \impl (A_0 \lor P(c))) =
\I(\tup{\Q}\tup{y}(B_0(\tup{d}) \impl (A_1 \lor P(c))))$ for all
$\tup{d}$, $c$.  To see this, first note that shifting quantifiers
across $\lor$, and shifting universal quantifiers out of the
consequent of a conditional is always possible. Hence it suffices to
consider the case of $\exists$.  $\I(\exists y\, A_2)$ is
either $= 0$ or $= 1$.  In the former case, both sides equal
$\I(B_0(\tup{d}) \impl P(d))$, in the latter, both sides equal~$1$.
\end{proof}

In summary, we obtain the following characterization of axiomatizability of
prenex fragments of Gödel logics:

\begin{theorem}
  The prenex fragment of $\gdl{V}$ is axiomatizable if and only if $V$ is
  finite or uncountable.  The prenex fragments of any two $\gdl{V}$ where $V$
  is uncountable coincide.
\end{theorem}

%
%

\subsection{$\bot$-free fragments\label{sec:botfree}}

\def\nobot{{\not\bot}}

In the following we will denote the $\bot$-free fragment of $\gdl{V}$ with
$\gdl{V}^\nobot$.  $\gdl{V}^\nobot$ is the set of all $\gdl V$-valid formulas
which do not contain $\bot$ (and hence also no $\lnot$).  First we show that
the only candidates for r.e. fragments are the $\bot$-free fragments
of $\gdl V$ where $V$ is uncountable and either $0\in V^\infty$ or $0$ is
isolated~$V$.

\begin{lemma}\label{botfree-notra}
  If $\gdl{V}$ is not r.e., then  $\gdl{V}^\nobot$ is also not r.e.
\end{lemma}

Define $A^b$ as the formula obtained from $A$ by replacing all occurences of
$\bot$ with the new propositional variable $b$ (a 0-place predicate symbol).
Then define $A^*$ as
\[ 
    A^* = \big(\lAnd_{P\in A}\qa\bar x(b\limp P(\bar x))\big)\limp A^b 
\]
where $P\in A$ means that $P$ ranges over all predicate symbols occuring
in~$A$.  We will first prove a lemma relating $A^*$ and $A$:

\begin{lemma}\label{bot-valid}
  \[ \gdl{V} \entails A \quad\text{iff}\quad \gdl{V}^{\nobot} \entails A^* \]
\end{lemma}
\begin{proof}
  If: Replace $b$ by $\bot$.
  
  Only if: Suppose $\gdl{V}^\nobot\nentails A^*$. Thus, there is an
  interpretation $\I_0$ such that $\I_0(A^*) < 1$. By
  Proposition~\ref{fact:loewenheim} and
  Lemma~\ref{lemma:interpretation-cut-off}, there is an interpretation $\I$
  such that $\I(A^b) < 1$ and $\I((\lAnd_{P\in A}\qa\bar x(b\limp P(\bar
  x))))=1$.  Because of the latter, for every atomic subformula~$B$ of $A$,
  $\I(B) \ge \I(b) = v$.  Define $\I'(B)$ for atomic subformulas $B$ of $A$ by
  \[
   \I'(B) = \begin{cases} 0 & \I(B) \le v\\
                       \I(B)& \I(B) > v
             \end{cases}
  \]
  (and arbitrary for other atomic formulas).  It is easily seen by induction
  that $\I'(B) = \I(B)$ if $\I(B) > v$, and if $\I(B) = v$, then $\I'(B) = v$
  or $=0$.  In particular, $\I'(A^b) < 1$.  But, of course, $\I'(b) = \I'(\bot)
  = 0$, and hence $\I'(A^b) = \I'(A)$.



\end{proof}

\begin{proof}[Proof of Lemma~\ref{botfree-notra}]
  If $\gdl{V}^\nobot$ were recursively enumerable, then by
  Lemma~\ref{bot-valid}, $\gdl{V}$ would also be recursively
  enumerable. 
\end{proof}

Thus, by Theorem~\ref{thm:count-nonax}, we only have two candidates for
axiomatizable $\bot$-free fragments: both truth-value sets have a non-empty
perfect kernel~$P$, and in the one case $0\in P$ and in the other $0\not\in P$
but $0$ is isolated. The prototypical Gödel sets for these cases are
$V_1 = [0,1]$ and $V_2 = \{0\} \cup [1/2,1]$. We will show that the $\bot$-free
fragments of these two logics coincide, thus in fact proving that there is
only one axiomatizable $\bot$-free fragment.
\fxnote{This should be done also for entailment!}
\begin{lemma}\label{embedding-half}
  Let $V_1 = [0,1]$ and $V_2 = \{0\} \cup [1/2,1]$.
  The $\bot$-free fragments of $\gdl{V_1}$ and $\gdl{V_2}$ coincide,
    i.e.\ 
    \[ \gdl{V_1}^\nobot\entails A\quad\text{iff}\quad
       \gdl{V_2}^\nobot\entails A
    \]
\end{lemma}
\begin{proof}
  Only if: obvious, since a counter-example in $V_2$ actually also is a
  counter-example in $V_1$.
  
  If: Suppose that $\gdl{V_1}^\nobot \nentails A$, i.e., there is an
  $\I_1$ such that $\I_1(A)<1$. Define $\I_2$ for all atomic
  subformulas~$B$ of $A$ by $\I_2(B) = 1/2(1+\I_1(B))$. By 
  Lemma~\ref{lem:induced-interpretation} and the remark following it
  we see that the definition of $\I_2$ extends to all formulas.
\end{proof}

\begin{theorem}
  The $\bot$-free fragment of $\gdl{V}$ is recursively axiomatizable
  if and only if $V$ is finite or uncountable and either $0$ belongs to
  $V^\infty$ or is isolated.  The $\bot$-free fragment of any two
  such $V$ coincide.
\end{theorem}
\fxnote{We should give the axiomatization}
\begin{proof}
  From Lemma~\ref{botfree-notra}, Lemma~\ref{embedding-half} and
  Theorem~\ref{thm:gs-id} for the
  uncountable case. The finite case is obvious as the additional axioms
  $\axfinn$ do not contain $\bot$.
\end{proof}

%
%
\subsection{$\forall$-free fragments}\label{sec:efrag}

In the following we will denote the $\exists$-fragment of $\gdl{V}$ with
$\gdl{V}^\exists$. It is the set of all formulas valid in $\gdl V$ which do not
contain~$\forall$.

First we show, as in the case of the $\bot$-free fragment, that the
only candidates for axiomatizable fragments are the two uncountable
ones, $0\in P$ and $0$ isolated.   We will do this by showing that the
formulas used to reduce validity in the other cases to Trachtenbrodt's Theorem
are validity-equivalent to $\forall$-free formulas.  

\begin{lemma}\label{lm:qshift}
  If $A(x)$ and $B$ are $\forall$-free, then
  \[ \entails \qa x\, A(x)\limp B \quad\text{iff}\quad \entails
  \qe x(A(x)\limp B)\]
\end{lemma}
\begin{proof}
  If: This is a valid quantifier shift
  rule. 
  
  Only if: Suppose that $\nentails\qe x(A(x)\limp B)$, i.e., there is an
  interpretation~$\I$ such that $\I(\qe x(A(x)\limp B)<1$. But this implies
  that
  \begin{equation}
    \label{eq:foo}
    \qa u\in U \quad \I(A(u))>\I(B).
  \end{equation}
  Now define $\I'(Q)$ for atomic subformulas~$Q$ of $A$ by 
  \begin{equation*}
    \I'(Q) = \begin{cases}
      \I(Q) & \text{if\ } \I(Q) \le \I(B) \\
      1     & \text{if\ } \I(Q) > \I(B).
    \end{cases}
  \end{equation*}
  Then (i) If $C$ is $\forall$-free and $\I(C) > \I(B)$, then $\I'(C) = 1$,
  and if $\I(C) \le \I(B)$, then $\I'(C) = \I(C)$; and (ii) $\I'(\qa x\,
  A(x)) = 1$
  
  (i) For atomic $C$ this is the definition of $\I'$. The cases for $\land$,
  $\lor$, and $\limp$ are trivial. Now let $C \equiv \qe x\, D(x)$. If $\I(\qe
  x\, D(x)) > \I(B)$, then for some $u \in U^\I$, $\I(D(u)) > \I(B)$. By
  induction hypthesis, $\I'(D(u)) = 1$ and hence $\I'(\qe x\, D(x)) = 1$.
  Otherwise, $\I(\qe x\, D(x)) \le \I(B)$, in which case $\I'(D(u)) =
  \I'(D(u))$ for all $u$.  (ii) By~(\ref{eq:foo}), for all $u\in U$, $\I(A(u))
  > \I(B)$, hence, by (i), $\I'(A(u)) = 1$.
 
  By (i) and~(ii) we have that $\I'(\qa x\, A(x)) = 1$ and
  $\I'(B) = \I(B) < 1$, thus $\I'(\qa x\, A(x) \limp B) < 1$, i.e., $\nentails
  \qa x\, A(x) \limp B$.
\end{proof}

Note that in the preceding Lemma we can replace the prefix of $A(x)$
by a string of universal quantifiers and the same proof will work.

\begin{lemma}\label{exists-notra}
  If $\gdl{V}$ is not recursively enumerable, then also $\gdl{V}^\exists$.
\end{lemma}

\begin{proof}
  It is sufficient to show that Formula~\ref{fm:countable} for
  $A^g$ as given on page~\pageref{fm:countable} and
  Formula~\ref{fm:0notiso} for $A^h$ as given on 
  page~\pageref{fm:0notiso} are validity-equivalent to $\forall$-free formulas.

  If we only consider the quantifier structure of these formulas and
  apply valid quantifier shifting rules, including the shifting rule
  for crisp formulas given in Lemma~\ref{lm:crisp}, we obtain in both
  cases formulas which are of the form
  \[ \qa \bar x A(\bar x) \limp B \]
  where $A(\bar x)$ and $B$ are $\forall$-free. By to Lemma~\ref{lm:qshift} we
  see that both formulas are validity equivalent to $\forall$-free formulas.
\end{proof}

As for the $\bot$-free fragments, it turns out that the two prototypical
examples of Gödel sets create the same $\exists$-fragment:

\begin{lemma}\label{embhalf2}
  Let $V_1 =[0,1]$ and $V_2 =\{0\}\cup[1/2,1]$.  The $\exists$-fragments of
  $\gdl{V_1}$ and $\gdl{V_2}$ coincide, i.e.\ 
  \[ \gdl{V_1}^\exists \entails A \quad\text{iff}\quad
  \gdl{V_2}^\exists \entails A
  \]
\end{lemma}
\begin{proof}
  Only if: obvious, since a counter-example in $V_2$ actually also is a
  counter-example in $V_1$.
  
  If: Suppose that $\gdl{V_1}^\exists \nentails A$, i.e., there is an $\I_1$
  such that $\I_1(A)<1$. Define $\I_2$ for all atomic subformulas~$B$ of $A$
  by $\I_2(B) = 1/2(1+ \I_1(B))$ if $\I_1(B) > 0$ and $= 0$ if $\I_1(Q)=0$.
  By Lemma~\ref{lem:induced-interpretation} and the remark following it
  we see that the definition of $\I_2$ extends to all formulas.
\end{proof}

\fxnote{This should be done for entailment if possible}
\begin{theorem}
  The $\exists$-fragment of $\gdl{V}$ is r.e. if and only if $V$ is finite or
  uncountable and either $0$ belongs to $V^\infty$ or is isolated.  The
  $\exists$-fragment of any two such $V$ coincide.
\end{theorem}
\fxnote{We should give a comment on the axiomatization!}
\begin{proof}
  From Lemma~\ref{exists-notra}, Lemma~\ref{embhalf2} and
  Theorem~\ref{thm:gs-id} for the uncountable case. The finite case is obvious
  as the additional axioms $\axfinn$ do not contain universal quantifiers.
\end{proof}

\section{Conclusion}

In the preceding sections, we have given a complete characterization of the
r.e. and non-r.e. first-order Gödel logics.  Our main result is that there are
two distinct r.e. infinite-valued G\"odel logics, viz., $\GR$ and $\Go$.  What
we have not done, however, is investigate how many \emph{non}-r.e. Gödel
logics there are.  It is known that there are continuum-many different
propositional consequence relations and continuum-many different propositional
quantified Gödel logics \cite{BaazVeith99LC}.  In forthcoming work \cite{bgp},
it is shown that there are only countably many first-order Gödel logics.
Although this result goes some way to clarifying the situation, a criterion of
identity of Gödel logics using some topological property of the underlying
truth value set is a desideratum.  We have only given
(Lemma~\ref{lem:g-embed}) a sufficient condition: if there is a continuous
bijection between $V$ and $V'$, then $\gdl{V} = \gdl{V'}$.  But this condition
is not necessary: any pair of non-isomorphic uncountable Gödel sets with $0$
contained in the perfect kernel provides a quick counterexample (as any two
such sets determine~$\GR$ as their logic).  Such a topological
characterization of first-order infinite valued Gödel logics could then be
used to obtain a more fine-grained analysis of the complexity of the non-r.e.
Gödel logics.  As noted already, these also differ in the degree to which they
are non r.e. \cite{Hajek05}.

Another avenue for future research would be to carry out the characterization
offered here for extensions of the language.  Candidates for such
extensions are the addition of the projection modalities 
($\bigtriangleup a = 0$ if $a=1$ and $=1$ if $a <1$), of the
globalization operator of \cite{TT86}, or of the involutive negation 
($\sim a = 1 - a$).  
It is known that $\GR$ with the addition of these operators is still
axiomatizable.  The presence of the projection modality, in particular,
disturbs many of the nice features we have been able to exploit in this paper,
for instance, in the presence of $\bigtriangleup$ the crucial
Lemma~\ref{lemma:interpretation-cut-off} and
Proposition~\ref{prop:consequences} no longer hold.  Thus, not all of our
results go through for the extended language and new methods will have to be
developed.

\bibliographystyle{alpha}
\bibliography{bpz}

\end{document}